
\input amstex
\documentstyle{amsppt}
\pagewidth{6in}
\NoRunningHeads
\NoBlackBoxes
\loadbold
\topmatter
\title A family of singular oscillatory integral operators \\
and failure of weak amenability
\endtitle
\author Michael Cowling, Brian Dorofaeff, Andreas Seeger, and
James Wright
\endauthor
\thanks
Research for this paper was supported in part by the Australian
Research Council and the National Science Foundation.
\endthanks
\address
Michael Cowling, School of Mathematics, University of New South
Wales, Sydney NSW 2052, Australia
\endaddress
\email m.cowling\@maths.unsw.edu.au\endemail
\address
Brian Dorofaeff, Department of Mathematics, Sydney Grammar School,
Darlinghurst, NSW 2010, Australia
\endaddress
\email bdd\@sydgram.nsw.edu.au\endemail
\address
Andreas Seeger, Department of Mathematics, University of
Wisconsin-Madison, Madison, WI 53706, USA
\endaddress
\email seeger\@math.wisc.edu\endemail
\address
James Wright, Department of Mathematics and Statistics,
University of Edinburgh, King's Building, Mayfield Road, Edinburgh EH3 9JZ, U.K.
\endaddress
\email wright\@maths.ed.ac.uk\endemail
\endtopmatter

\define\tu{{\tilde u}}
\define\tv{{\tilde v}}

\define\tw{{\tilde w}}
\define\tz{{\tilde z}}

\define\R{\Bbb{R}}
\define\N{\Bbb{N}}
\define\Z{\Bbb{Z}}

\define\group#1{\text{\rm #1}}
\define\Lie#1{\frak{#1}}
\define\Heis{\group{H}}
\define\SLtwoR{\group{SL}(2,\R)}

\define\Symp{\group{Sp}}
\define\SO{\group{SO}}
\define\SU{\group{SU}}

\define\ad{\operatorname{ad}}
\define\Ad{\operatorname{Ad}}

\define\transpose{^T}
\define\betsqb{(1 + b^2/4)}
\define\VN{V\!N}

\define\MA{M\!A}
\define\MzA{M_0A}
\define\Proj{\operatorname{Proj}}

\redefine\vec#1{#1}
\define\End{\operatorname{End}}
\define\spaN{\operatorname{span}}

\define\mfr#1{{\frak{#1}}}

\define\calD{\Cal{D}}

\define\calH{\Cal{H}}

\define\calM{\Cal{M}}
\define\calO{\Cal{O}}

\define\calS{\Cal{S}}

\define\veczeta{\zeta}
\define\piez{\pi_{\eta,\veczeta}}
\define\chiez{\chi_{\eta,\veczeta}}
\define\Hez{\calH_{\eta,\veczeta}}

\document

\define\sgn{\operatorname{sign}}
\define\supp{\operatorname{supp}}
\define\inn#1#2{\langle#1,#2\rangle}

\define\lc{\lesssim}


\define\ep{\epsilon}

\define\ka{\kappa}
            
\define\la{\lambda}             
\define\sig{\sigma}

\define\om{\omega}

\define\fS{{\frak S}}

\define\fg{{\frak g}}


\define\fr{{\frak r}}
\define\fs{{\frak s}}

\define\bbR{{\Bbb R}}

\define\bbT{{\Bbb T}}

\define\bbZ{{\Bbb Z}}

\define\cB{{\Cal B}}

\define\cD{{\Cal D}}
\define\cE{{\Cal E}}

\define\cH{{\Cal H}}

\define\cK{{\Cal K}}

\define\cM{{\Cal M}}

\define\cO{{\Cal O}}
\define\cP{{\Cal P}}
\define\cQ{{\Cal Q}}
\define\cR{{\Cal R}}
\define\cS{{\Cal S}}
\define\cT{{\Cal T}}
\define\cU{{\Cal U}}
\define\cV{{\Cal V}}
\define\cW{{\Cal W}}

\define\cZ{{\Cal Z}}



\define\set#1{\{{#1}\}}
\define\emph#1{{{\it #1}}}
\define\norm#1{\|{#1}\|}
\define\pviint{{\text{p.v.}\iint}}
\define\pvint{{\text{p.v.}\int}}

\head{\bf 1.  Weak amenability and $\boldsymbol \Lambda
\boldkey(\boldkey{G}\boldkey)$}\endhead

Let $G$ be a locally compact Hausdorff topological group, equipped
with a left-invariant Haar measure, written $dx$ or $dy$ in
integrals. We write $L^p(G)$ for the usual Lebesgue space of
(equivalence classes of) functions on~$G$. In this section, the
symbol $\lambda$ will denote the left regular representation of
$G$ on $L^2(G)$, and $f*g$ will denote the convolution of
functions $f$ and $g$ on $G$.

\subhead
{\bf 1.1. The Fourier algebra and pointwise multipliers}
\endsubhead

A matrix coefficient of the left regular representation is a
function of the form
$$
x \mapsto \inn{\lambda(x)h}{k}
  =\int h(x^{-1}y) \overline {k(y)} \,dy
$$
where $h$ and $k$ lie in $L^2(G)$. The {\it Fourier algebra\/} of
$G$, denoted by $A(G)$, is defined to be the Banach space of all
these, that is,
$$
A(G) = \set{ \inn{\lambda(\cdot) h}{k}: h,k\in L^2(G) }
$$
(which is actually a linear space), equipped with the norm
$$
\|\varphi\|_A = \inf\set{ \|h\|_2 \|k\|_{2} : \varphi =
\inn{\lambda(\cdot) h}{k}}.
$$
The infimum is in fact attained, see \cite{12}.
If $G$ is abelian, a function in $A(G)$ is the Fourier transform
of a function in $L^1(\widehat G)$, where $\widehat G$ is the dual
group of~$G$.

All functions in $A(G)$ are continuous and vanish at infinity. The
Fourier algebra forms a commutative Banach algebra under pointwise
operations, with Gel'fand spectrum $G$. It has a unit (the
function~$1$) if and only if $G$ is compact.  For proofs of these
results and for much more information about the Fourier algebra,
see the original article by
 Eymard \cite{12} or the book by Pier \cite{29}.

The \emph{group von Neumann algebra} $\VN(G)$ is defined to be the
set of all bounded linear operators on $L^2(G)$ commuting with
right translations.  Suppose that $f \in L^1(G)$.  We associate to
$f$ the left convolution operator $\lambda[f]$ on $L^2(G)$,
defined by
$$
\lambda[f] h(x) = \int_G \lambda(y) h(x) \, f(y) \,dy = f * h(x).
$$
This operator lies in $\VN(G)$. The function $f$ also gives rise
to an element of the topological dual space $A(G)^*$ of $A(G)$, by
integration: one defines $L_f$ in $A(G)^*$ to be the linear
functional $\varphi \mapsto \int_G \varphi(x) \, f(x) \, dx$.  The
association between the operator $\lambda[f]$ and the linear
functional~$L_f$ extends to identify the group von Neumann
algebra~$\VN(G)$ with the dual space $A(G)^*$.
More precisely, for any $F$ in $A(G)^*$,
there exists a unique $F'$ in $\VN(G)$ such that
$$
\inn{F'(h)}{k} = F(\inn{ \lambda (\cdot)h }{k})
\qquad\forall h,k\in L^2(G).
$$
The mapping $F\mapsto F'$ is an isometric isomorphism; it also
carries the weak-star topology of $A(G)^*$ to the ultraweak
topology of $\VN(G)$. The set $\set{L_f : f\in L^1(G)}$ is
weak-star dense in $A(G)^*$ and the set $\set{\lambda[f]:f\in
L^1(G)}$ is ultraweakly dense in~$\VN(G)$.  The correspondence
between $F$ and $F'$ is the unique continuous extension of the map
$\lambda[f] \mapsto L_f$. For proofs of these facts, see \cite{12}
or \cite{29}.

On a Lie group~$G$, $\calD(G) \subset A(G)$, where $\calD$ denotes
the space of compactly supported smooth functions. We may think of
elements of~$A(G)^*$ as distributions on~$G$, and of elements of
$\VN(G)$ as convolutions by these distributions.

We shall need the notion of a completely bounded operator on a von
Neumann algebra.  Suppose that $\calM$ is a von Neumann algebra
and $T:{\calM}\to{\calM}$ is a continuous linear operator. Let
$\calM^{(n)}$ be the algebra of $n\times n$ matrices with entries
in $\calM$ and let $I_n$ be the $n\times n$ identity matrix.
Define the extension $T\otimes I_n$ to $\calM^{(n)}$ by $(T\otimes
I_n F)_{ij}= T(F_{ij})$. Then $T$ is said to be \emph{completely
bounded} if
$$
c_{_T} = \sup_n\| T\otimes I_n\|<\infty .
$$
We write $\|T\|_{cb}$ for the completely bounded operator norm~$c_{_T}$. 
Much more about completely bounded operators may
be found in \cite{28}.

We define $\MA(G)$, the \emph{space of (pointwise) multipliers} of
$A(G)$, to be the set of all continuous functions $\varphi$ on $G$
such that the pointwise product $\varphi\psi$ lies in $A(G)$ for
all $\psi$ in $A(G)$.  A multiplier $\varphi\in \MA(G)$ may be
identified with the multiplication operator $m_\varphi$ on $A(G)$
given by $m_\varphi: \psi\mapsto\varphi\psi$, and we equip
$\MA(G)$ with the corresponding operator norm.

We also define $\MzA(G)$, the \emph{space of completely bounded
multipliers of $A(G)$}, also called Herz--Schur multipliers (see,
e.g., \cite{4}), to be the set of all continuous
functions $\varphi$ on $G$ such that the adjoint operator
$m^*_\varphi$ is completely bounded as an operator on $\VN(G)$. We
define $\|\varphi\|_{\MzA(G)}$ to be the completely bounded
operator norm~$\|m^*_\varphi\|_{cb}$.  This space is smaller than
$\MA(G)$, and the norm is larger than the $\MA(G)$-norm. For
further information about these spaces, see the articles by
 Cowling \cite{4} and
De~Canni\`ere and Haagerup \cite{8}; in particular it is shown in \cite{8}
 that
$\|m\|_{M_0A(G)}=\sup_H \|m\otimes 1_H\|_{MA(G\times H)}$ where
the supremum is taken over all locally compact groups $H$.

Both $\MA(G)$ and $\MzA(G)$ form Banach algebras under pointwise
multiplication. We have the inclusions $A(G)\subseteq B(G)
\subseteq \MzA(G) \subseteq \MA(G)$ where $B(G)$ is the
Fourier--Stieltjes algebra consisting of matrix coefficients of
unitary representations. If the group~$G$ is {\it amenable}, i.e.,
there exists a left invariant mean on $L^\infty(G)$, then both of
these algebras coincide with the Fourier--Stieltjes algebra $B(G)$;
in fact the equality $B(G)=MA(G)$ is a characterization of
amenability, see \cite{26}. In general, these inclusions are proper;
in fact, specific examples of functions in $M_0A(G)$ arise as
matrix coefficients of uniformly bounded representations which
need not be equivalent to unitary ones (see  \cite{25}, \cite{30}).

\subhead{\bf 1.2.
Approximate units}
\endsubhead

Let $L$ be a positive real number. Then $A(G)$ is said to have an
\emph{approximate unit bounded by $L$} if there exists a directed
set $I$ and a net $\set{\varphi_i: i\in I}$ of functions in $A(G)$
such that
$$
\lim_{i\in I}\|\psi-\varphi_i\psi\|_{A}=0
\qquad\forall\psi\in A(G)
\tag 1.2.1
$$
and
$$
\|\varphi_i\|_{A}\le L\qquad\forall i\in I.
\tag 1.2.2
$$
It is known that $A(G)$ has an approximate unit bounded by a
positive real number $L$ if and only if $A(G)$ has an approximate
unit bounded by $1$; this is one of the many equivalent conditions
for $G$ to be \emph{amenable} (Leptin \cite{24}, see also
Herz \cite{17}).  When $G$ is
amenable, the existence of the approximate unit implies that
$$
\norm {\varphi }_{A} = \norm {\varphi }_{\MzA}
= \norm {\varphi }_{\MA}
\qquad\forall \varphi \in A(G).
$$
For more information about amenability, see \cite{29}.

One may weaken the existence criterion on the approximate unit as
follows. Given a positive real number $L$, we say that $A(G)$ has
an \emph{$L$-completely bounded approximate unit}, if there exists
a net $\set{\varphi_i: i\in I}$ of functions in $A(G)$ such that
(1.2.1) holds and
$$
\|\varphi_i\|_{\MzA}\le L\qquad\forall i\in I
\tag 1.2.3
$$
We define the number $\Lambda(G)$ to be the infimum of all the
numbers $L$ for which there exists an $L$-completely bounded
approximate unit on $A(G)$, with the convention that
$\Lambda(G)=\infty$ if no such approximate unit exists. The
group~$G$ is said to be \emph{weakly amenable} if
$\Lambda(G)<\infty$.

Finally we say that $A(G)$ has an \emph{$L$-multiplier bounded
approximate unit}, if there is a net $\set{\varphi_i: i\in I}$ of
functions in $A(G)$ such that (1.2.1) holds and
$$
\|\varphi_i\|_{\MA}\le L\qquad\forall i\in I.
$$
A multiplier bounded
approximate unit  is simply an
$L$-multiplier bounded
approximate unit, for some $L<\infty$.

Clearly $\Lambda(G)\in [1,\infty]$, because
$\|\cdot\|_\infty\leq\|\cdot \|_{\MzA(G)}$, but in every known
case, $\Lambda(G)$ is an extended integer. Much of what is known
about $\Lambda(G)$ for locally compact groups is summarized in
the following list.  For details see the articles
by Haagerup \cite{13}, \cite{14}, Cowling \cite{4}, \cite{5},
 De~Canni\`ere and  Haagerup \cite{8},
Cowling and Haagerup \cite{6},
Lemvig Hansen \cite{23},
 Bo\D zejko and  Picardello \cite{1}, Dorofaeff \cite{9}, \cite{10}.

\proclaim{1.2.1}
Suppose that $G$, $G_1$, and $G_2$ are locally compact groups.

 (i)
If $G_1$ is isomorphic to $G_2$, then $\Lambda(G_1)=\Lambda(G_2)$.

 (ii)
If $K$ is a compact normal subgroup of $G$, then
$\Lambda(G)=\Lambda(G/K)$.

 (iii)
If $G_1$ is a closed subgroup of $G_2$, then $\Lambda(G_1) \le
\Lambda (G_2)$, with equality if $G_2/G_1$ admits a finite $G_2$
invariant measure.

 (iv)
If $G$ is the direct product group $G_1\times G_2$, then
$\Lambda(G)=\Lambda(G_1) \, \Lambda(G_2)$.

 (v)
If $G$ is discrete and $Z$ is a central subgroup of~$G$, then
$\Lambda(G)\le\Lambda(G/Z)$.

 (vi)
If $G$ is amenable, then $\Lambda(G)=1$.

 (vii)
If $G$ is a free group, then $\Lambda(G)=1$.

 (viii)
If $G$ is an amalgamated product $G=*_A G_i$, where each $G_i$ is
an amenable locally compact group, and $A$ is a compact open
subgroup of all $G_i$, then $\Lambda(G)=1$.

 (ix)
If $G$ is locally isomorphic to $\SO(1,n)$ or to $\SU(1,n)$, then
$\Lambda(G)=1$.

 (x)
If $G$ is locally isomorphic to $\Symp(1,n)$, then
$\Lambda(G)=2n-1$.

 (xi)
If $G$ is locally isomorphic to $F_{4(-20)}$, then
$\Lambda(G)=21$.

 (xii)
If $G$ is a simple Lie group of real rank at least two, then
$\Lambda(G)=\infty$.
\endproclaim

For generalizations of these ideas to von Neumann algebras, see
Haagerup \cite{13}, \cite{14}, Cowling and Haagerup \cite{6}
and for generalizations to ergodic systems
and dynamical systems, see Cowling and Zimmer \cite{7} and
Jolissaint \cite{19}.  These ideas are
loosely related to Property (T) and the Haagerup Property,
which are investigated in detail in the books
by Zimmer \cite{36},  by de la Harpe and Valette \cite{16} and by
Ch\'erix,  Cowling,  Jolissaint,  Julg and  Valette \cite{3}.

We shall make use of the following results, without further
reference.

\proclaim{1.2.2}
Suppose that $H$ is a closed subgroup of the locally compact
group~$G$, that $T$ is a distribution on~$H$, and that $\varphi$
is a function on~$G$. Then:

 (i)
if $\varphi\in A(G)$, then $\varphi\bigr|_H\in A(H)$
 and $\|\varphi\bigr|_H\|_{A(H)} \leq \|\varphi\|_{A(G)}$

(ii)
if $T \in A(H)^*$, then, considered as a distribution on~$G$,
$T \in A(G)^*$ and $\|T\|_{A(G)^*} = \|T\|_{A(H)^*}$

 (iii)
if $\varphi\in \MzA(G)$, then $\varphi\bigr|_H\in \MzA(H)$ and
$\|\varphi\bigr|_H\|_{\MzA(H)}\le\|\varphi\|_{\MzA(G)}$

 (iv)
if $\varphi\in \MA(G)$, then $\varphi\bigr|_H\in \MA(H)$ and
$\|\varphi\bigr|_H\|_{\MA(H)}\le\|\varphi\|_{\MA(G)}$.
\endproclaim

\noindent See \cite{17, Thm.~1} and
\cite{8, Prop.~1.12} for the proofs.

\proclaim{1.2.3}
If $\set{\varphi_i: i\in I}$ is an $L$-completely bounded
approximate unit on $A(G)$, then $\varphi_i\to 1$ uniformly on
compact subsets of $G$. Conversely, if there exists a net
$\set{\varphi_i: i\in I}$ of $A(G)$-functions such that $
\|\varphi_i\|_{\MzA(G)}\le L$ and $\varphi_i\to 1$ uniformly on
compact sets, then there exists an $L$-completely bounded
approximate unit of compactly supported $A(G)$-functions, $\set{
\tilde\varphi_j: j\in J}$ say. If $G$ is a Lie group, then we may
also assume that $\tilde\varphi_j\in \calD(G)$ for all $j$ in~$J$.

This result also holds when ``$L$-completely bounded'' is replaced
by ``$L$-multiplier bounded''.
\endproclaim
\noindent
For the proof, see \cite{6, Prop.~1.1}.

\proclaim{1.2.4}
 Let $K$ be a compact normal subgroup of the locally compact group $G$.

(i) Let $m\in MA(G)$ and define for $\widetilde m(gK)=\int_K m(gk) dk$
(where $dk$ is normalized Haar measure). Then $\widetilde m\in MA(G/K)$ with
$\|\widetilde m\|_{MA(G/K)}\le \|m\|_{MA(G)}$.

(ii) The statement (i) remains true for $MA(G)$ replaced with $MA_0(G)$;
moreover the space $MA_0(G/K)$ may be isometrically
identified with the subspace of functions in $MA_0(G)$
which are constants on the cosets of $K$ in $G$.
Furthermore $\Lambda(G/K)=\Lambda(G)$.
\endproclaim

\noindent  (i) is immediate. For (ii)
see \cite{6, Prop.~1.3} (one uses the definition
\cite{6, (0.3)} to verify the nontrivial part of (ii)).

\proclaim{1.2.5}
Suppose that $G=SK$ is a (set) decomposition of $G$ as a product
of an amenable closed subgroup $S$ and a compact subgroup $K$, and
that $\nu$ is normalized Haar measure on $K$. Suppose further that
$\tilde A(G)$ is one of $A(G)$ or $\MzA(G)$ or $\MA(G)$. Then for
any $\varphi\in \tilde A(G)$ the average $\dot\varphi$, defined by
$$
\dot\varphi(x)=\int_{K\times K}\varphi(k x k') \,d\nu(k)\,d\nu(k'),
$$
belongs to $\tilde A(G)$. Further, $\|\dot\varphi\|_{\tilde
A(G)}\le \|\varphi\|_{\tilde A(G)}$.
\endproclaim
\noindent For the proof, see \cite{6, Prop.~1.6}.
The point of the lemma is that, by
averaging, we may assume that any given approximate unit of
$A(G)$-functions bounded in the $\tilde A(G)$-norm is
$K$-biinvariant, with the same bound. The above lemma also holds
if we choose compactly supported smooth functions, and these
properties are preserved by averaging.

\subhead
{\bf What lies ahead}
\endsubhead

For a connected noncompact simple Lie group $G$ with finite center
and real rank at least two, the invariant $\Lambda(G)$ takes the
value infinity. This result was proved by Haagerup \cite{14}.
His proof involves investigating certain semidirect products,
namely $\SLtwoR\ltimes\R^2$ and $\SLtwoR\ltimes \Heis^1$, where
$\Heis^1$ is the Heisenberg group of dimension three. He shows that
these semidirect products do not admit multiplier bounded
approximate units, and hence deduces that $\Lambda$ is infinite
for both the semidirect products and then, by structure theory,
for any noncompact simple Lie group $G$ with finite center and
real rank at least two. These semidirect products are the smallest
members of two families of semidirect products, for which it turns
out to be interesting to calculate $\Lambda$ (see Section~8). The
first family is formed with the action of the unique irreducible
representation of $\SLtwoR$ on $\R^n$.  It was shown by Dorofaeff
\cite{9} that all these groups have infinite~$\Lambda$; this was used to show that
that the original hypothesis  of finite center in Haagerup's proof of
1.2.1 (xii) is redundant (\cite{10}). The
second family is where $\SLtwoR$ acts on the Heisenberg group
$\Heis^n$ of dimension~$2n+1$ by fixing the center and operating
on the vector space $\R^{2n}$ by the unique irreducible
representation of dimension~$2n$.

We consider this family of semidirect products and show they do
not admit multiplier bounded approximate units; in particular
$\Lambda(\SLtwoR\ltimes \Heis^n)=\infty$. Given this and earlier
results, and some structure theory, it is now possible to compute
$\Lambda(G)$ for any real algebraic Lie group~$G$, or indeed for
any Lie group~$G$ whose Levi factor has finite center.

\proclaim{Main Theorem}
Let $G$ be a real Lie group with Lie algebra~$\fg$, and let
$\fs\oplus \fr$ be the Levi decomposition of~$\fg$, where $\fr$ is
the maximal solvable ideal of $\fg$ and $\fs$ is a semisimple
summand, and let $\fs_1\oplus\dots\oplus\fs_m$ be the
decomposition of $\fs$ as a sum of simple ideals. Let $S$ be a
maximal analytic semisimple subgroup of $G$ corresponding to
$\fs$, and let $S_i$ be the subgroup associated to $\fs_i$, where
$i=1,\dots,m$. Suppose that $S$ has finite center.

Then $G$ is weakly amenable if and only if one of the following
two conditions is satisfied for each $i=1,\dots,m$:

\noindent Either

(*) $S_i$ is compact

\noindent or

(**) $S_i$ is noncompact, of real rank $1$, and the action of
$\fs_i$ on $\fr$ is trivial, i.e., $[\fs_i,\fr]=0$.

\medskip

If for every $i\in\{1,\dots,m\}$, either (*) or (**) is satisfied
then $\Lambda(G)=\prod_{i=1}^m\Lambda(S_i)$ and $\Lambda(G)$ can
be computed by consulting the list (1.2.1).

If for at least one $i\in\{1,\dots,m\}$ neither (*) nor (**)
holds, then $A(G)$ does not admit any multiplier bounded
approximate unit.
\endproclaim

{\it Structure of the paper.}
The main part of this paper (Sections~2--7) is devoted to the
proof that the Fourier algebra of
$\SLtwoR\ltimes \Heis^n$ does not admit multiplier bounded
 approximate units, and consequently we have
 $\Lambda(\SLtwoR\ltimes \Heis^n)=\infty$. Using a
modification of Haagerup's approach for the case  $n=1$ \cite{14}, one
can reduce matters to the estimation of a singular oscillatory
integral operator; this reduction is described in Section~2. The
estimation of the integral operator, which is rather nontrivial,
is carried out in Sections~3--7. In Section~8 we consider general
Lie groups under the assumption that the Levi part has finite center.
Here we use facts from  the structure theory of  Lie groups to show that
if  for at least one $i\in\{1,\dots,m\}$ neither condition (*) nor condition (**) in the
Theorem holds, then $G$ does not admit multiplier bounded approximate units.
This  will be combined  with previously known results to complete the  proof of the main theorem.

\head{\bf 2. A family of semidirect products}\endhead

Fix a positive integer $n$.  Throughout this chapter we shall
consider the group
$$
G_n=\SLtwoR\ltimes \Heis^n,
$$
where $\SLtwoR$ acts on the
Heisenberg group  $\Heis^n$ by the
unique irreducible representation of dimension $2n$,   fixing  the
center. We shall reduce the proof that $\Lambda(G_n)=\infty$ to
the estimation of a family of singular oscillatory integral
operators. The four subsequent sections will then be dedicated to
estimating these operators.

\subhead{\bf 2.1. The action of $ \SLtwoR$ on the Heisenberg
group}
\endsubhead

Recall that $\Heis^n$ is a Lie group whose
underlying manifold is $\R^{2n}\times\R$.  The group
multiplication may be given by the formula
$$
(\vec u,t)(\vec u',t')
 =(\vec u + \vec u', t + t' + \vec u^T B \vec u'),
$$
where the symplectic matrix $B$ is defined by
$$
B _{ij}= \cases (-1)^j & \text{if $i+j=2n+1$} \\
  0   & \text{otherwise.}
\endcases
$$
We shall write $\{\vec e_1, \dots, \vec e_{2n}\}$ for the standard
basis of $\R^{2n}$.

We shall now describe the action of $\SLtwoR$ on $\R^{2n}$ by the
irreducible representation~$\pi_{2n}$ of dimension $2n$ which is unique up to isomorphism
(see, e.g.,  \cite{22, p.~107}). For $j =
1, \dots, 2n$, let $$\alpha_j=\binom {2n-1}{j-1}^{1/2}.$$ We
identify $\R^{2n}$ with the space $\Bbb P_{2n}$ of homogeneous
polynomials in two variables of degree $2n-1$ by associating
$(u_1,\dots,u_{2n})$ with the polynomial
$$
 P: (x,y) \mapsto \sum_{j=1}^{2n}\alpha_j u_j x^{2n-j} y^{j-1},
\tag 2.1.1
$$
and define the action of $A$ in $\SLtwoR$ by
$$
\pi_{2n}(A)P(x,y)= P((x,y)A)= P(ax+cy,bx+dy)
\qquad\forall (x,y)\in\R^2,
$$
where $A=\pmatrix a&b\\c&d\endpmatrix$ (see \cite{20}). If $P$
is as in (2.1.1), then a computation shows that
$$
\pi_{2n}(A)P(x,y)
=\sum_{i=1}^{2n}\bigl[Z(A)u\bigr]_i \alpha_i x^{2n-i}y^{i-1},
$$
where the $2n\times 2n$ matrix $Z(A)$ is given by
$$
\bigl(Z(A)\bigr)_{ij}=
\sum_{l=0}^{2n} \binom{j-1}{l} \binom{2n-j}{2n-i-l}
\alpha_i^{-1} \alpha_j\,
     a^{2n-i-l} \, b^l \, c^{i+l-j} \, d^{j-l-1}
\tag 2.1.2
$$
(see \cite{9}).   Here we use the standard
convention that $\binom k l = 0$ if $l$ is negative or $l>k$.

In order to extend the action on $\R^{2n}$ to an action on
$\Heis^n$ we need to show that the action on $\R^{2n}$ is
symplectic.

\proclaim{Lemma 2.1.1}
The map $Z$ is a symplectic action on $\R^{2n}$, i.e.,
$$
Z(A)\transpose \, B \, Z(A)= B
\tag 2.1.3
$$
for each $A\in\SLtwoR$. Define $\Bar Z(A):\Heis^n\to \Heis^n$ by
$$
\Bar Z(A) (\vec u,t)=(Z(A)\vec u,t);
$$
then $\Bar Z(A)$ is an automorphism of $\Heis^n$ and $\Bar Z$ is
an action of $\SLtwoR$ on $\Heis^n$.
\endproclaim
\demo{\bf Proof}
Recall that $\alpha_j=\binom {2n-1}{j-1}^{1/2}$. {}From (2.1.2)
and our choice of $\alpha$ one checks that
$$
Z(A^T)=Z(A)^T.
$$
Observe also that
$$
B=Z(J) \qquad\text{where}\qquad J = \pmatrix 0&1\\-1&0
\endpmatrix.
$$
For any $A$ in $\SLtwoR$, a direct matrix calculation shows
$A\transpose JA=J$ and so
$$
Z(A)\transpose \, B Z(A) = Z(A)\transpose Z(J)Z(A)
=Z(A\transpose J A)=Z(J)
$$
and therefore (2.1.3) holds. The fact that $\Bar Z(A)$ is an
automorphism of $\Heis^n$ follows immediately from (2.1.3); hence
$\Bar Z$ is an action on $\Heis^n$. \qed
\enddemo

We may now describe the semidirect product group $G_n$. As a
manifold, this is $\SLtwoR\times\R^{2n}\times\R$. The product in
$G_n$ is defined by
$$
(A,\vec u,t)(A',\vec u',t')=(AA',\vec u+Z(A)\vec u',t+t'+\vec u^T
BZ(A)\vec u')
$$
and the inverse is given by
$$
(A,\vec u,t)^{-1}=(A^{-1},-Z(A^{-1})\vec u,-t),
\tag 2.1.4
$$
for all $(A,\vec u,t)$ and $(A',\vec u',t')$ in $G_n$. The closed
subgroups $\{(I,\vec u, t): \vec u \in \R^{2n},\ t \in \R\}$
(where $I$ is the identity of $\SLtwoR$) and $\{(A,0,0) :
A\in\SLtwoR\}$ may be identified with $\Heis^n$ and $\SLtwoR$.
Given $(A',\vec u,t)$ and $(A,0,0)$ in $G_n$, it follows that
$$
(A,0,0) (A',\vec u,t) (A,0,0)^{-1} = (AA'A^{-1}, Z(A)\vec u,t),
$$
which shows that $\Heis^n$ is normalized by $\SLtwoR$.

There are several important subgroups and elements of $G_n$ which
we now identify.  We denote by $K$ the compact subgroup
$\SO(2,\R)$ of $\SLtwoR$, considered as a subgroup of $G_n$.  For
$b$ in $\R$, we define
$$
k_b^\pm = \pm \beta(b)^{-1}
\pmatrix b/2& 1 \\ - 1 & b/2 \endpmatrix,
\qquad
n_b =
\pmatrix 1 & b \\ 0 & 1 \endpmatrix,
\qquad\text{and}\qquad
h_b =\pmatrix 0 & \beta(b) \\ -\beta(b)^{-1} & 0
\endpmatrix,
$$
where $\beta(b) = \betsqb^{1/2}$. Then $k_b^\pm \in K$. We write
$N$ for the nilpotent subgroup $\set{ n_b : b\in\R }$. For future
purposes, we observe the following lemma.

\proclaim{Lemma 2.1.2}
For all $b$ in $\R$, we have
$$
\align
k^+_b \, n_b \, k^-_b &= n_{-b},
\\
n_{b/2} \, k^+_b \, n_{b/2} &= h_b.
\tag 2.1.5
\endalign
$$
Further,
$$
Z(h_b)(u_n \vec e_n + u_{n+1} \vec e_{n+1}) = (-1)^n
(\beta(b)^{-1} u_n \vec e_{n+1} - \beta(b) \, u_{n+1} \vec e_n).
$$
Finally,
$$
\bigl(Z(n_b)\bigr)_{ij}=\cases \alpha_i^{-1} \alpha_j \binom
{j-1}{j-i} b^{j-i} &\quad\text{if } j>i
\\
1 &\quad\text{if }  j=i
\\
0 &\quad\text{if }  j<i
\endcases
\tag 2.1.6
$$
and, in particular,
$Z(n_b)_{n,n+1} = nb.
$
\endproclaim
\demo{\bf Proof}
These are all straightforward computations which will be omitted.
\qed\enddemo

\subhead{\bf 2.2. Two nilpotent subgroups}\endsubhead

We write $G$ for $G_n$, and $H$ for the subgroup of $G$ of all
elements of the form $(n_b, \vec u, t)$, where $b \in \R$ and
$\vec u \in \R^{2n}$. Let $V_k$ denote the subspace~$\spaN\{
\vec e_1, \dots, \vec e_k \}$ of $\R^{2n}$ (when $k=1,\dots, 2n$).
Since $N$ is a subgroup of $\SLtwoR$ and the matrix $Z(n_b)$ is
upper triangular for all $b$ in $\R$, this subspace is invariant
under all the maps $Z(n_b)$, and the subset of $G$ of all elements
of the form $(n_b, \vec v, t)$, where $b \in \R$ and $\vec v \in
V_k$, is a subgroup of $H$. We write $H_0$ for the subgroup of $G$
obtained in this way when $k = n+1$.

We need to understand the behavior of the restrictions of
$K$-bi-invariant functions on $G$ to $H$.
It
follows from formula (2.1.4) that
$$
   (k_b^+,0,0) (n_b,\vec u,t) (k_b^-,0,0)
 = (n_{-b}, Z(k_b^+)\vec u,t).
\tag 2.2.1
$$
We define the diffeomorphism $\Omega: H \to H$ by the formula
$$
\Omega(n_b,\vec u,t) = (n_{-b}, Z(k_b^+)\vec u,t).
\tag 2.2.2
$$

\proclaim{Lemma 2.2.1}
If $\varphi \in \calD(G)$ and $\varphi$ is $K$-bi-invariant, then
$\varphi|_H\circ\Omega = \varphi|_H$.
\endproclaim
\demo{\bf Proof} Since $\varphi$ is $K$-bi-invariant, we have
$$\varphi(n_b,\vec u,t)
 = \varphi\bigl((k_b^+,0,0) (n_b,\vec u,t) (k_b^-,0,0)\bigr)$$
for all $(n_b,\vec u,t)$ in $H$
and the assertion follows from formulae (2.2.1) and (2.2.2).
\qed
\enddemo

\subhead{\bf 2.3. Some distributions on $\boldkey H_{\boldkey 0}$}
\endsubhead

We will define a family of distributions on $H_0$, using two
iterated principal value integrals.  To clarify the sense in which
these are to be interpreted, and because it will be useful later,
we first discuss certain principal value integrals on $\R^2$.  For
Schwartz functions $\psi\in \calS(\R^2)$, let
$$
\aligned D(\psi) &= \pviint \frac{\psi(s_1,s_2)}{s_2^2 - s_1^2}
\,\,ds_1\,ds_2
\\
&= \lim_{\epsilon \to 0+}
        \lim_{\delta \to 0+}
\int_{|s_2|>\epsilon}  \frac{1}{2 s_2}
\Bigl(
\int_{|s_2+s_1|>\delta} \frac{\psi(s_1,s_2)}{s_2+s_1} ds_1\,+\,
\int_{|s_2-s_1|>\delta} \frac{\psi(s_1,s_2)}{s_2-s_1} ds_1\Bigr) ds_2.
\endaligned
\tag 2.3.1
$$
It is routine to show that $D$ is a tempered distribution. We
shall also need a modification $\widetilde D$ defined by
$$
\widetilde D(\psi)=D(\widetilde \psi) \text{ where } \widetilde
\psi(y_1,y_2)=\psi(y_2,y_1).
\tag 2.3.2
$$
The distributions $D$ and $\widetilde D$ satisfy
$$
D(\psi)+\widetilde D(\psi)=\pi^2\psi(0,0)
\tag 2.3.3
$$
for all Schwartz functions; this fact was used by Haagerup \cite{14}
and called the failure of Fubini's theorem,
since it can be rewritten in the form
$$
\pviint\frac{\psi(s_1,s_2)}{s_2^2-s_1^2} \bigl(ds_1ds_2-ds_2ds_1\bigr)=\pi^2
\psi(0,0).
$$
The verification of formula (2.3.3) can be found in \cite{9};
it relies on a Fourier transform
calculation and the fact that $D( e^{-i\inn{\cdot}{\tau}})$ is
equal to $\pi^2$ if $\tau_1^2>\tau_2^2$ and to $0$ if
$\tau_1^2<\tau_2^2$.

For fixed $b\in \R$, define $Q_b:\R^2\to H_0\subset G$ by
$$
Q_b(s_1,s_2)=(n_b,Z(n_{b/2})(s_1e_n+s_2 \beta(b)^{-1}e_{n+1}),0)
\tag 2.3.4
$$
where, as before, $\beta(b) = \betsqb^{1/2}$.  For a test function
$\phi\in\calD(G)$, let $Q_b^*\phi$ be the pullback of $\phi$ to
$\R^2$ defined (as usual) by $Q_b^*\phi(s_1,s_2)=
\phi(Q_b(s_1,s_2))$. We now define the distribution $D_R$ on $G$,
for all $R$ in $\R^+$, by the formula
$$
D_R(\phi)=\int_{-R}^R
D(Q_b^*\phi) \frac{db}{\beta(b)}
\tag 2.3.5
$$
We may view $D_R$ as a distribution on $H$ or on $G$, with support
in $H_0$, if we wish.

\proclaim{Lemma 2.3.1}
Suppose that $\varphi \in \calD(G)$ and $\varphi$ is
$K$-bi-invariant.  Then
$$
D_R(\varphi) = \frac{\pi^2}2 \int_{-R}^R
  \frac{\varphi(n_b,\vec 0,0)}{\betsqb^{1/2}} \,\,db.
\tag 2.3.6
$$
In particular, if $\{\varphi_n\}_{n\in \N}$ is a sequence of
$K$-bi-invariant $\calD(G)$-functions, and $\varphi_n \to 1$
uniformly on compact subsets of $G$ as $n \to \infty$, then
$$
\lim_{n\to \infty} D_R(\varphi_n) = 2 \pi^2
 \sinh^{-1}(R/2).
\tag 2.3.7
$$
Both formulae remain valid if $D_R$ is considered as a
distribution on $H$ or $H_0$ and applied to restrictions of
$K$-bi-invariant functions to $H$ or $H_0$.
\endproclaim
\demo{\bf Proof}
Recall from Lemma 2.2.1 that if $\varphi \in \calD(G)$ and
$\varphi$ is $K$-bi-invariant, then $\varphi|_H =
\varphi|_H\circ\Omega$. Now we compute for arbitrary $\phi\in
\cD(G)$
$$
\align
\phi\circ \Omega(Q_b(s_1,s_2))
&=
\phi(n_{-b}, Z(k_b^+)Z(n_{b/2})(s_1e_n+s_2\beta(b)^{-1}e_{n+1}),0)
\\
&=
\phi(n_{-b}, Z(n_{-b/2})Z(h_b)(s_1e_n+s_2\beta(b)^{-1}e_{n+1}),0)
\\
&=
\phi(n_{-b}, Z(n_{-b/2})
(-1)^n(\beta(b)^{-1}s_1e_{n+1}-s_2e_{n})
,0) .
\endalign
$$
Here we have used the definition of $\Omega$ and $Q_b$, and the
relation $k_b^+n_{b/2}=n_{-b/2} h_b$ (which follows from formula
(2.1.5)). Since $D$ is even on $\R^2$, it follows that
$$
D(Q_b^*(\phi\circ\Omega) )=\widetilde D(Q_{-b}^*\phi),
$$
and therefore, since $\beta$ is even,
$$
D_R(\phi\circ\Omega)= \int_{-R}^R \widetilde D(Q_{b}^*\phi)
\beta(b)^{-1} db .
\tag 2.3.8
$$

Now we assume that $\varphi$ is $K$-bi-invariant and use (2.3.3).
Then
$$
\align
D_R(\varphi)&= \frac 12
\int_{-R}^R
\bigl( D(Q_{b}^*\varphi)
+\widetilde D(Q_{b}^*\varphi)\bigr) \beta(b)^{-1} db
\\
&=\frac{\pi^2}2
\int_{-R}^R
Q_{b}^*\varphi(0,0)
 \beta(b)^{-1} db
\\&=\frac{\pi^2}2
\int_{-R}^R \varphi(n_b,0,0) \beta(b)^{-1} db ,
\endalign
$$
as required.

The formula (2.3.7) follows by passing to the limit and evaluating
the integral.

The last assertion follows from our computation, since $\Omega$
maps the subset of $G$ (or of $H$ or $H_0$) consisting of all
$(n_b, s_1e_n+s_2e_{n+1},0)$ into itself. \qed
\enddemo

\subhead{\bf 2.4. Failure of weak amenability}\endsubhead

We are now in a position to reduce the question of the weak
amenability of $G$ to a question of boundedness of the operators
$\la[D_R]$ of convolution with~$D_R$.

\proclaim{Proposition 2.4.1}
Suppose that $\lambda[D_R]$ lies in $\VN(H_0)$, and that
$\norm{\lambda[D_R]} = o(\log R)$ as $R\to\infty$. Then $G$ is not
weakly amenable, i.e., $\Lambda(G) = \infty$, and further, there
does not exist a multiplier bounded approximate unit on $G$.
\endproclaim

\demo{\bf Proof}
If $G$ were weakly amenable, then there would exist $L$ in $[1,
\infty)$ and a sequence $\set{ \varphi _n : n \in \N }$ of
$\calD(G)$-functions such that $\norm{ \varphi_n }_{\MzA} \leq L$
for all $n$ in $\N$ and $\varphi_n \to 1$, uniformly on compact
subsets of $G$, as $n \to \infty$.  By averaging if necessary, we
could suppose that all the functions $\varphi_n$ were
$K$-bi-invariant; see (1.2.5).  \emph{A fortiori}, for some $L$ in~$\R^+$, there
would be a sequence $\set{ \varphi_n: n\in \N}$ of
$K$-bi-invariant $\calD(G)$-functions satisfying the conditions
$\norm{ \varphi_n }_{\MA} \leq L$ and $\varphi_n \to 1$ as $n \to
\infty$.  The same would be true if there existed a multiplier
bounded approximate unit on~$G$.

Consider the sequence $\set{ D_R(\varphi_n|_{H_0}) : n \in \N }$.
Since $H_0$ is amenable, $A(H_0)$ has an approximate unit, whence
$$
\norm{ \varphi_n|_{H_0} }_{A} =  \norm{ \varphi_n|_{H_0} }_{\MA}
\leq \norm{ \varphi_n }_{\MA} \leq L.
\tag 2.4.1
$$
Thus
$$
\big| D_R(\varphi_n|_{H_0}) \big| \leq \norm{ \lambda[D_R] }_{\VN}
\norm{ \varphi_n |_{H_0} }_{A} \leq L \norm{\lambda[ D_R] }_{\VN}  = o(\log
R).
$$
However, by (2.3.7)
$$
\lim_{n \to \infty} \big| D_R(\varphi_n|_{H_0}) \big| = 2 \pi^2
\log \bigl(\tfrac R 2+\sqrt{\tfrac{R^2}4 +1}\bigr).
$$
The last two formulae are contradictory, so the original
hypothesis of the weak amenability of $G$ must be incorrect.
\qed\enddemo

Most of this paper is  dedicated to verifying the hypothesis of Proposition 2.4.1; 
more precisely, we shall obtain  the estimate
$$
\norm{\lambda[D_R] }_{\VN(H_0)}=O(\log\log R)
\quad \text{ as
$R\to \infty$.}
\tag 2.4.2
$$
To do this, we will use Fourier analysis on $H_0$ to study the
distributions $D_R$ when acting on $A(H_0)$. The first stage in
this process is to find a family of unitary representations $\set{
\piez : \eta \in \R, \veczeta \in V_n}$ of $H_0$; we then describe
the Plancherel formula for this group. It is a consequence of the
Plancherel formula that $\norm{ \lambda[D_R]}_{\VN}$ is equal to
the supremum of the operator norms $\norm{ \piez [D_R] }$ as
$\eta$ and $\veczeta$ vary. We shall then identify the
operators~$\piez[D_R]$ as singular oscillatory integral operators,
which will be estimated in Sections~3--7.

\subhead {\bf 2.5. Representations of the group $\boldkey H_{\boldkey 0}$}
\endsubhead

To simplify notation, from now on we write $(b, u, t)$ instead of
$(n_b, u, t)$, and $P(b)$ instead of $Z(n_b)$, see (2.1.6). Then
the group law may be rewritten in the form
$$
(b,\vec u,t)(b',\vec u',t')=(b+b',\vec u+P(b) \vec u',t+t'+\vec u^T
BP(b) \vec u')
$$
and
$$
(b,\vec u,t)^{-1}=(-b,-P(-b)\vec u,-t),
$$
for all $(b,\vec u,t)$ and $(b',\vec u',t')$ in $H_0$. {}From
formula (2.1.3), it follows that $P(-b)^T B P(-b) = B$, so
$$
\aligned
     (b, \vec u,t)^{-1} (b', \vec u', 0)
 & = (b' - b, P(-b)(\vec u' - \vec u), -t -\vec u^T B \vec u') \\
 & = (b' - b, P(-b)(\vec u' - \vec u), -t + (-1)^n(u_nu_{n+1}' - u_{n+1}u_n')).
\endaligned
\tag 2.5.1
$$

It is easy to see that the subgroup $H_1$ of $H_0$, given by
$$
H_1 = \set{ (0, \vec w, s) : \vec w \in V_n,\ s \in \R } ,
$$
is normal in $H_0$ and abelian. Let $\fS$ be the subset $\set{ (c,v
\vec e_{n+1},0) \in H_0 : c,v \in \R }$ of $H_0$.  As a set, we
may identify $\fS$ with $\R^2$. Any element $h$
of $H_0$ may be expressed uniquely in the form $\sig h_1$, where $\sig
\in \fS$ and $h_1 \in H_1$.  Indeed, if $c$, $s$, $t$, and $v$ are
in $\R$, while $\vec w \in V_n$ and $\vec u \in V_{n+1}$, then
$$
(c, v\vec e_{n+1}, 0)(0, \vec w, s) = (c, v\vec e_{n+1} + P(c)
\vec w, s + v \vec e_{n+1}^TB P(c) \vec w) ,
$$
so
$$
\gathered
(c, v\vec e_{n+1}, 0)(0, \vec w, s) = (b, \vec u, t)
\qquad\text{if and only if}
 \\
  c = b , \quad
  v = u_{n+1} , \quad
 \vec w = P(-b)\Proj_{V_n} \vec u , \quad\text{and}\quad
  s = t - (-1)^n u_{n+1} u_{n},
\endgathered
\tag 2.5.2
$$
where $\Proj_V $ denotes the standard orthogonal projection onto
the subspace $V$ of $\R^{2n}$. As a consequence, we also note the
integration formula
$$
\int_{H_0}F(y) dy=\int_{H_1}\int_{\fS} F(\sigma z) d\sigma dz.
\tag 2.5.3
$$

We define the characters $\chiez $ of $H_1$ by the formula
$$
\chiez  (0,\vec w, s)=\exp (i(-1)^n\eta s + i\inn{\veczeta}{ \vec w}),
\tag 2.5.4
$$
where $\eta \in \R$ and $\veczeta \in V_n^*$, and induce the
character $\chi_{-\eta,-\zeta}$ from $H_1$ to $H_0$. The induced
representation $\pi_{\eta,\zeta}$ acts on the Hilbert space $\Hez$
of all complex-valued functions $\xi$ on $H_0$ such that
$$
\xi((b,\vec u,t)(0,\vec w,s))=\chiez (0,\vec w,s)\xi(b,\vec u,t)
\qquad\forall (0, \vec w, s) \in H_1
\quad\forall (b,\vec u,t) \in H_0 ,
$$
and
$$
\Bigl(\int |\xi(c,v\vec e_{n+1},0)|^2 \,dc \,dv\Bigr)^{1/2}<\infty.
$$
We equip this space with the norm equal to the left hand side of
this inequality. As $H_0 = \fS\,H_1$, each function in $\Hez$ is
determined by its restriction to~$\fS$, and so this really is a
norm on $\Hez$, modulo the usual issues of identification of
functions which differ on null sets. Clearly $\calH_{\eta,\zeta}$
can be identified with $L^2(\fS)$.

The action of the unitary representation $\piez $ on a
function $\xi$ in $\Hez$ is defined by the formula
$$
\piez (b,u,t)\xi (b', u', t')=\xi\bigl((b,u,t)^{-1}(b',u',t')\bigr).
$$
In particular, using formulae (2.5.1) and (2.5.2)
and
we see that
$$
\align
\piez (b, \vec u,t)\xi(c, v\vec e_{n+1},0)
& =\xi\bigl((b, \vec u,t)^{-1}(c, v\vec e_{n+1}, 0)\bigr) \\
& =\xi\bigl(c-b, P(-b)(v\vec e_{n+1} - \vec u), -t +(-1)^n u_n v \bigr)\\
& =\xi\bigl((c-b, (v-u_{n+1})\vec e_{n+1}, 0)(0,\vec w,s)\bigr),
\endalign
$$
where $(0,\vec w,s)$ in $H_1$ is defined by
$$\align
\vec w &= P(b-c)\Proj_{V_n} P(-b)(v\vec e_{n+1} - \vec u) \\
       &= P(b-c)\bigl[
            P(-b)(v\vec e_{n+1} - \vec u) - (v-u_{n+1})\vec e_{n+1} \bigr] \\
       &=   P(-c)(v\vec e_{n+1} - \vec u) + P(b-c)(u_{n+1}-v)\vec
                e_{n+1},
\endalign
$$
and since  $P(b)_{n,n+1} = n b$ by Lemma 2.1.2,
$$\align
s  &= -t + (-1)^nu_nv - (-1)^n (v - u_{n+1})(P(-b)_{n,n+1}(v - u_{n+1}) - u_n)
\\
  &= -t + (-1)^n\bigl(nb (v - u_{n+1})^2 + u_n (2v  - u_{n+1}) \bigr).
\endalign
$$
In conclusion,
$$
\align
&\piez (b, \vec u,t)\xi(c, v\vec e_{n+1},0)
\tag 2.5.5
\\
& =\xi(c-b, (v-u_{n+1})\vec e_{n+1}, 0)
\\
& \qquad\times
        \exp\bigl(i\eta\bigl[(-1)^{n+1}t +nb(v-u_{n+1})^2 +u_n(2v-u_{n+1})
        \bigr]\bigr)
\\
& \qquad\times
        \exp \bigl(i
\inn{\veczeta}{P(-c)(v \vec e_{n+1} - \vec u)
                +P(b-c)(u_{n+1} - v)\vec e_{n+1}} \bigr) .
\endalign
$$
The elements of $\fS$ act by
translations (here we think of~$\fS$ as~$\R^2$), combined with
multiplications, while the action of the elements of~$H_1$ is as
follows:
$$\gather
\piez (0, \vec w,t)\xi(c, v\vec e_{n+1},0)
\\
 =\xi(c, v\vec e_{n+1}, 0)
    \exp(i\eta[(-1)^{n+1}t + 2w_nv] - i\inn{\veczeta}{  P(-c)\vec w}).
\endgather
$$

Finally we extend the representation $\piez$ to functions $f$ in
$L^1(H_0)$. For each $\eta$ in $\R$ and $\veczeta$ in $V_n^{*}$,
we associate an operator $\piez [f]$ on $L^2(\fS)$ in the usual
way by the formula
$$
\piez [f] \xi (\sigma) = \int_{H_0} f(x) \,\piez (x) \xi(\sigma) \,dx .
\tag 2.5.6
$$
This formula extends by continuity to define a Fourier transform
of certain distributions on $H_0$.

%
 \subhead{\bf 2.6 A Plancherel formula}\endsubhead

In what follows we shall write $\chi$ for $\chi_{\eta,\zeta}$ and
$\cH_\chi$ for $\cH_{\eta,\zeta}$; we also denote by $d\chi$ the
measure $(2\pi)^{-n-1} d\eta \,d\zeta$ on the dual space $\widehat
H_1$.

For $\Xi$ in $\cD(H_0)$ and $\chi$ in $\widehat H_1$, define the
function $\Xi_\chi$ on $H_0$ by
$$
\Xi_\chi(x)=\int_{H_1}\Xi(xz)\overline{\chi}(z) dz.
\tag 2.6.1
$$
We note that

\proclaim{Lemma 2.6.1}
For all $\Xi$ in $\cD(H_0)$, the function $\Xi_\chi$ belongs to the
Hilbert space $\calH_\chi$. Further
$$
\|\Xi\|_{L^2(H_0)}=\Bigl(\int_{\widehat H_1}
\|\Xi_\chi\|^2_{\cH_\chi}d\chi\Bigr)^{1/2}
$$
and the map $\Xi\mapsto\bigl(\chi\mapsto \Xi_\chi)$ extends to an
isometric bijection of $L^2(H_0)$ to $L^2(\widehat H_1,\cH)$.
\endproclaim
\demo{\bf Proof}
For $\Xi \in \cD(H_0)$ and $\chi\in \widehat H_1$, we compute:
$$\align
\Xi_\chi(xz')&=\int_{H_1}\Xi(xz'z)\overline{\chi}(z) dz
=\int_{H_1}\Xi(xz)\overline{\chi}({z'}^{-1}z) dz
\\&=\chi(z') \int_{H_1}\Xi(xz)\overline{\chi}(z) dz=
\chi(z')\Xi_\chi(x),
\endalign
$$
so that $\Xi_\chi$ has the required covariance property. Further
as $\sigma$ varies over $\fS$, the function $\Xi_\chi(\sigma)$
varies smoothly, and as a function on $\fS$ it has compact
support, contained in $\supp(\Xi)H_1\cap \fS$. Moreover by the
Plancherel theorem for $H_1$, and Fubini's theorem,
$$
\align
\int_{\widehat H_1}\|\Xi_\chi\|^2_{\cH_\chi} d\chi
&=
\int_{\widehat H_1}\int_\fS|\Xi_\chi(\sigma)|^2 d\sigma d\chi
=\int_\fS\int_{\widehat H_1}|\Xi_\chi(\sigma)|^2 d\chi d\sigma
\\
&=\int_\fS\int_{ H_1}|\Xi(\sigma z)|^2 dz d\sigma
=\int_{H_0}|\Xi(y)|^2 dy.
\endalign
$$
The extension to $L^2(H_0)$ is straightforward.\qed\enddemo

\proclaim{Lemma 2.6.2}
Suppose that $D$ is a distribution in $H_0$ and suppose that the
operator norm on $L^2(\cH_\chi)$ satisfies $\|\pi_\chi[D]\|\le A $
for all $\chi\in \widehat H_1$. Then $\lambda[D]$ is in $VN(H_0)$
and $\|\lambda[D]\|_{VN}\le A$.
\endproclaim
\demo{\bf Proof}
We shall assume that $D$ is given by integration against a
$\cD(H_0)$ function $k$; the general case follows by a
regularization argument. Now let $\Xi$ and $\Gamma$ be in
$L^2(H_0)$. Then
$$
\align
\inn{\la[k]\Xi}{\Gamma}
&=
\int_{H_0}\int_{H_0}k(x)\Xi(x^{-1}y)\overline\Gamma(y)dy dx
\\&=
\int_{H_0}\int_\fS\int_{H_1}k(x)\Xi(x^{-1}\sigma z)
\overline\Gamma(\sigma z)dz d\sigma dx
\\&=
\int_{H_0}\int_\fS\int_{\widehat H_1}k(x)\Xi_\chi(x^{-1}\sigma )
\overline{\Gamma_\chi(\sigma )}d\chi d\sigma dx
\\&=
\int_{\widehat H_1}\int_\fS
\int_{H_0}
k(x)\Xi_\chi(x^{-1}\sigma )
\overline{\Gamma_\chi(\sigma )}dx d\sigma d\chi
\\&=
\int_{\widehat H_1}
\int_\fS \pi_\chi[k]\Xi_\chi(\sigma )
\overline{\Gamma_\chi(\sigma )} d\sigma d\chi
\\&=
\int_{\widehat H_1} \inn{ \pi_\chi[k]\Xi_\chi}{
\Gamma_\chi}_{\cH_\chi} d\chi,
\endalign
$$
where we used (2.5.3), the Plancherel theorem on the abelian group
$H_1$, Fubini's theorem and the definitions of $\pi_\chi(f)$ and
$\cH_\chi$. {}From the hypothesis and the Cauchy--Schwarz
inequality it follows that
$$
\align
\big|\inn{\la[k]\Xi}{\Gamma}\big|
&\le \int_{\widehat H_1}| \inn{ \pi_\chi[k]\Xi_\chi}{
\Gamma_\chi}_{\cH_\chi}| d\chi
\\&\le
\int_{\widehat H_1} A
\| \Xi_\chi\|_{\cH_\chi}\|\Gamma_\chi\|_{\cH_\chi} d\chi
\\&\le A
\Bigl(
\int_{\widehat H_1}
\| \Xi_\chi\|_{\cH_\chi}^2
d\chi\Bigr)^{1/2}
\Bigl(
\int_{\widehat H_1}
\| \Gamma_\chi\|_{\cH_\chi}^2
d\chi\Bigr)^{1/2}
\\
&=A\|\Xi\|_{L^2(H_0)}
\|\Gamma\|_{L^2(H_0)},
\endalign
$$
by Lemma 2.6.1.  Taking the supremum over all $\Xi$ and $\Gamma$ with  norm $\le 1$
shows that $\|\la[k]\|_{VN}\le A$. \qed
\enddemo

%
\subhead{\bf 2.7. The oscillatory singular integral operators
$\boldsymbol{\pi}_{\boldsymbol\eta,\boldsymbol
\zeta}\boldkey[\boldkey D_{\boldkey R} \boldkey] $}
\endsubhead

We now compute the operator-valued Fourier transform of the
distributions $D_R$.  We change notation slightly, and for $\xi$
in $\Hez$, we write $\xi(c,v)$ instead of $\xi(c,v\vec
e_{n+1},0)$. We also set
$$
q(b)=n^{-1}\inn{ \veczeta}{P(b)\vec e_{n+1}},
\tag 2.7.1
$$
and write $\calM_q$ for the operator on $L^2(\fS)$ of pointwise
multiplication by the function $(c,v) \mapsto \exp(invq(-c))$.
Observe that
$$
\align
   q'(b)
&=\frac{\partial}{\partial b}\,n^{-1}\,
\inn{\veczeta}
{        \sum^{n+1}_{i=1} P(b)_{i,n+1}\,\vec e_i}\\
&= n^{-1}\,\frac{\partial}{\partial b}\,
\inn{\veczeta}
{        \sum^{n+1}_{i=1} \binom{n}{ n+1-i} \binom{2n-1}{ n}^{1/2}
                          \binom{2n-1}{i-1}^{-1/2} b^{n+1-i}\,\vec e_i}\\
&=\inn{\veczeta}{ \sum^n_{i=1} \binom{n-1}{ n-i} \binom{2n-1}{ n-1}^{1/2}
        \binom{2n-1}{i-1}^{-1/2} b^{n-i}\,\vec e_i}\\
&=\inn{\veczeta}{ \sum^n_{i=1} P(b)_{in}\, \vec e_i}
                = \inn{\veczeta}{  P(b) \vec e_n}.
\endalign
$$
Then, rewriting formula (2.5.5)
we have shown that
$$
\align
& \piez (b, \vec u,t)\xi (c, v) \\
&=\xi(c-b, v-u_{n+1})
        \exp\bigl(i\eta\bigl[(-1)^{n+1}t +nb(v-u_{n+1})^2
                +u_n(2v-u_{n+1}) \bigr]\bigr)\\
& \qquad\times
        \exp \bigl(i\inn{\veczeta}{P(-c)(v \vec e_{n+1} - \vec u )
                + P(b-c)(u_{n+1} - v)\vec e_{n+1} } \bigr).
\endalign
$$
Thus
$$
\align
& \piez \bigl(b, P(b/2)(u_n \vec e_n + u_{n+1} \vec e_{n+1}) ,0\bigr)\xi (c, v) \\
&=\xi(c-b, v-u_{n+1})
        \exp\bigl(i\eta\bigl[nb(v-u_{n+1})^2 +(u_n + u_{n+1}nb/2)(2v-u_{n+1}\bigr)
                                \bigr]\bigr)\\
&\qquad \times\exp \bigl(i\inn{\veczeta}
{ P(-c)(v \vec e_{n+1} - P(b/2)
        (u_n \vec e_n + u_{n+1} \vec e_{n+1}) )
                + P(b-c)(u_{n+1} - v)\vec e_{n+1}
} \bigr) \\
&=\xi(c-b, v-u_{n+1})
        \exp\bigl(i\eta\bigl[ n b(v^2 - v u_{n+1} + u_{n+1}^2/2) +u_n (2v-u_{n+1})
                                \bigr]\bigr)\\
& \qquad\times
        \exp \bigl(i\bigl[n v q(-c)
        -  u_n q'(b/2-c) - n u_{n+1} q(b/2-c)
        +  n(u_{n+1} - v) q(b-c) \bigr] \bigr) ,
\endalign
$$
and so
$$
\align
&\qquad \piez [D_R]\xi (c, v) \\
&=\int_{-R}^R \pviint
        \piez \bigl(b, P(b/2)(u_n \vec e_n +\tfrac{ u_{n+1}}{\beta(b)}
 \vec e_{n+1}) ,0\bigr)\xi(c, v)
        \,\frac{1} {u^2_{n+1}-u^2_n} \, du_n \, du_{n+1} \,\frac{ db}{\beta(b)} \\
&=\int_{-R}^R \pviint
        \piez \bigl(b, P(b/2)(u_n \vec e_n + u_{n+1} \vec e_{n+1}) ,0\bigr)\xi(c, v)
        \,\frac{1} {\beta(b)^2 \, u^2_{n+1}-u^2_n} \,\, du_n \, du_{n+1} \, db \\
&=\int_{-R}^R \pviint \xi(c-b, v-u_{n+1})
        \exp\bigl(i\eta\bigl[ n b(v^2 - v u_{n+1} + u_{n+1}^2/2) +u_n (2v-u_{n+1})
                                \bigr]\bigr) \\
& \qquad\times
        \exp \bigl(i\bigl[n v q(-c)
        +  n(u_{n+1} - v) q(b-c)
        -  u_n q'(b/2-c) - n u_{n+1} q(b/2-c) \bigr] \bigr) \\
& \qquad\times
        \,\frac{1} {\beta(b)^2 \, u^2_{n+1}-u^2_n} \,\,du_n \,du_{n+1}
        \,db,
\endalign
$$
and consequently
$$
\align
&\qquad \cM_q^{-1}\piez [D_R]\cM_q\xi (c, v) \\
&=  \int_{-R}^R \pviint \xi(c-b, v-u_{n+1})
        \exp \bigl(-  iu_n q'(b/2-c) + i\eta u_n (2v-u_{n+1})\bigr) \\
& \qquad\times
        \exp\bigl(i\bigl[ n\eta b(v^2 - v u_{n+1} + u_{n+1}^2/2) - n u_{n+1} q(b/2-c)
                                \bigr]\bigr)
        \frac{1} {\beta(b)^2 \, u^2_{n+1}-u^2_n} \,\, du_n \, du_{n+1} \, db .
\endalign
$$
We can calculate the innermost integral exactly: indeed
$$
\align
\pvint\exp (i\lambda z)\,\frac{dz}{w^2-z^2}
& =\frac{1}{2w}\pvint\exp (i\lambda z) \Bigl[\frac{1}{z+w}
        -\frac{1}{z-w}\Bigr] \,dz
\\
& =- \frac{\exp (i\lambda w)-\exp (-i\lambda w)}{2w} \pvint
\exp (i\lambda z)\,\frac{dz}{z}
\\
& =\frac{\pi\sgn(\lambda)\sin (\lambda w)}{w} =\frac{\pi\sin (|\lambda| w)}{w}.
\endalign
$$
We deduce that
$$
\aligned
&\cM_q^{-1}\piez [D_R]\cM_q\xi (c, v) \\
&= \pi  \int_{-R}^R \pvint   \xi(c-b, v-u_{n+1})
        \sin \bigl( \beta(b)u_{n+1}|(\eta (2v-u_{n+1}) - q'(b/2-c))|\bigr) \\
& \qquad\times
        \exp\bigl(in\bigl[ \eta b(v^2 - v u_{n+1} + u_{n+1}^2/2) - u_{n+1} q(b/2-c)
                                \bigr]\bigr)
        \frac{1} {\beta(b) \, u_{n+1} } \,\,  du_{n+1} \, db
\endaligned
\tag 2.7.2
$$
Since the sine term vanishes when $u_{n+1}$ vanishes, the
principal value of the inner integral is the usual integral.

\subhead{2.8  Equivalent formulation of the oscillatory integrals}
\endsubhead

In (2.7.2), we make the change of variables $y_1 = c - b$, $y_2 =
v - u_{n+1}$, $x_1 = c$, and $x_2 = v$, and set $p(t)=2q(-t/2)$
(so that $p'(t) = -q'(t/2)$). Then
$$
\aligned
&\quad
\Cal M_q\pi_{\eta,\zeta}[D_R] \Cal M_q^{-1}\! f (x_1,x_2) \\
&= \pi  \iint_{\R^2}   f(y_1, y_2)
     \sin \bigl( \beta(x_1 - y_1) (x_2 - y_2)|\eta (x_2 + y_2) + p'(x_1+y_1)|\bigr) \\
& \times
        \exp\bigl(\frac{in}{2}\bigl[ \eta (x_1 - y_1)(x_2^2 + y_2^2)
                        - (x_2 - y_2) p(x_1 + y_1) \bigr]\bigr)
  \frac{\chi_{[-R,R]}(x_1-y_1)}{\beta(x_1 - y_1) \, (x_2 - y_2)}
\,\,dy_1 \,dy_2.
\endaligned
\tag 2.8.1
$$

Now if $\eta\neq 0$, one can conjugate with a dilation in the
second variable by a factor of~$|\eta|^{1/2}$, to reduce to the
case where $\eta = \pm 1$. Further, changing the sign of $\eta$
and of the polynomial $p$ has the effect of changing the kernel to
its complex conjugate, and a kernel operator is bounded on $L^2$
if and only if the operator with conjugate kernel is.  In short,
to establish the uniform boundedness of the operators
$\pi_{\eta,\zeta}[D_R]$, as $(\eta, \veczeta)$ varies
over~$\R\times V_n$, we may suppose that $\eta$ is equal to $1$.

\head{\bf 3. The oscillatory integral}
\endhead

\demo{\bf Notation}
{}From now on, fix a positive integer $n$ and $\Gamma$ in $(1/2,
\infty)$. Let $p$ be a real polynomial of degree at most~$n$. An
{\it admissible\/} constant means a constant which depends only on
$n$ and $\Gamma$. We write $A\lc B$ if $A \le C B$ and $C$ is an
admissible constant in this sense. All ``constants'' $C$ below
will be admissible, and may vary from place to place.
\enddemo

Define the functions $\Psi:\R^2 \to \R$ and $\theta:\R^2 \to \R$
by the formulae
$$
\align
\Psi(x,y)&= (x_1-y_1)(x_2^2+y_2^2)-(x_2-y_2)p(x_1+y_1)
\tag 3.1\\
\theta(x,y)&=\beta(x_1-y_1)|x_2+y_2+p'(x_1+y_1)|(x_2-y_2) ;
\tag 3.2
\endalign
$$
further, recall that
$$
\beta(t) = (1+t^2/4)^{1/2}.
$$
Suppose $\min\{1, n/2\} \le |\gamma|\le\Gamma$ (the relevant value
of $\gamma$ will be $n/2$).  For $R>0$, we define the family of
singular oscillatory integral operators $\calO^R$ by
$$
\calO^R f(x)=\iint \frac{e^{i\gamma\Psi(x,y)}\sin\theta(x,y)}{
\beta(x_1-y_1) (x_2-y_2)} \chi_{[-R,R]}(x_1-y_1) f(y_1,y_2)
\, dy_1 \, dy_2
\tag 3.3
$$
for all $f$ in $C^\infty_0(\R^2)$. We shall see easily that
$\|\Cal O^R\|_{L^2\to L^2}=O(\log R)$ as $R\to \infty$ (this
follows from Lemma 3.3 below). However, we have the following
result.

\proclaim{Theorem 3.0}
Suppose that $\min\{ 1, n/2\}\le |\gamma|\le \Gamma$ and $R\ge
100$. Then the operator~$\calO^R$ extends to a bounded operator on
$L^2(\R^2)$, and
$$
\|\calO^R\|_{L^2(\R^2)\to L^2(\R^2)}
\leq C_{n,\Gamma}\log\log(10+R),
$$
where $C_{n,\Gamma}$ is admissible. If $n=1$, then this estimate
may be improved to $\|\calO^R\|=O(1)$.
\endproclaim

It is conceivable that the bound $\|\calO^R\|=O(1)$ holds in the
general case, but this has not been proved so far. For our
application, the assertion of the Theorem is (more than) enough.

\proclaim{Corollary 3.1}
(i) $\sup_{\eta,\zeta}\|\pi_{\eta,\zeta}[D_R]\|_{\cH_{\eta,\zeta}\to
\cH_{\eta,\zeta}}\lc \log(\log (10+R))$.

(ii) $\|\la[D_R]\|_{VN(H_0)}\lc \log\log(10+R)$.

(iii) The Fourier algebra of $\SLtwoR\ltimes\Heis^n$ does not
admit multiplier bounded approximate units.
\endproclaim

\demo{Proof}
By the results of Subsection~2.8, it suffices to prove (i) with
$\eta=1$, but then formula (2.8.1) shows that the statement is
implied by Theorem~3.0. The calculations in Section~2.7 together
with Lemma 2.6.2 show that (i) implies (ii), and (ii) implies
(iii), by Proposition 2.4.1.
\enddemo


\remark{Remarks}

(i) The assumption $\min\{ 1, n/2\}\le |\gamma|\le \Gamma$ in
Theorem 3.0 can be replaced by $ 1/2\le |\gamma|\le \Gamma$.
However the proof for the case where $|\gamma|=1/2$ and $n\ge 2$
turns out to be substantially more complicated. Fortunately this case
is irrelevant for our application.

(ii) There are many results concerning singular oscillatory
integral operators with kernels of the form $k(x-y)e^{i P(x,y)}$,
where $P$ is a polynomial.  If $k$ is a standard
Calder\'on--Zygmund kernel, the oscillatory variants are $L^p$
bounded ($1<p<\infty$), see Ricci and Stein \cite{31}. If $k$ is a {\it
multiparameter\/} Calder\'on--Zygmund kernel the technique in
\cite{31}, which uses induction on the degree of the polynomial,
no longer applies. In fact the $L^2$ boundedness may then hold or
fail depending on the properties of the polynomial $P$; see, e.g.,
\cite{2},  where a complete characterization of boundedness is
obtained for the special case where $P(x,y)=q(x-y)$ and $q$ is a
polynomial of two variables.

No theory for general polynomials is currently available.
Moreover, our operator is not included in the general class of
operators just discussed, because of the positivity of $\beta$.
Our proof of Theorem 3.0 relies on a subtle global cancellation
property of the distribution $D$ defined in (2.3.1) and the
noncommutativity of the convolution structure.

(iii)  It is instructive to examine the analogue of $\lambda[D_R]$ in the 
commutative  setting, where we 
identify $H_0$ as a set with $\R^{n+3}$, writing $(b,u,t)$
for~$(n_b,u,t)$, and 
replace the matrix $Z(b)$ by the identity throughout. Thus
define $q_b:\bbR^2\to \bbR^{n+3}$ by $q_b(s_1,s_2)=(b, s_1e_n+s_2\beta(b)^{-1} e_{n+1},0)$ and 
define a distribution by
$\frak D_R(\phi) =\int_{-R}^R D(q_b^*\phi)\beta(b)^{-1}db $. Denote by $\Cal C_R f$ the
convolution $\frak D_R*_E f$ on $\R^{n+3}$; here $*_E$ refers to the
standard commutative convolution in Euclidean space. The
operator $\Cal C_R$ is bounded on $L^2(\R^{n+3})$; however there
is a lower bound for the operator norm of the form $\|\Cal
C_R\|\ge c\log R$ as $R\to \infty$. This can be quickly seen by applying the 
partial Fourier transform $\Cal F_{n+2}$ in the $(u,t)$-variables. Indeed
for fixed $(\xi,\tau)\in \bbR^{n+1}\times \bbR$ let $\Cal C_R^{\xi,\tau}$ be the operator on 
$L^2(\Bbb R)$ of convolution with $\Cal F_{n+2}[\frak D_R] (\cdot, \xi,\tau)$; then the operator norm
of $\Cal C_R$ is equal to $\sup_{\xi,\tau}\|\Cal C_R^{\xi,\tau}\|_{L^2(\bbR)\to L^2(\bbR)}$.
A quick calculation using the formula for the Fourier transform of $D$ mentioned in \S2.3 
shows that $\|\Cal C_R^{\xi,\tau}\|\lc \log R$ and in particular
$$
\Cal C_R^{\xi,\tau} g(b)=\pi^2\int_{-R}^R g(b-c) \frac{dc}{\beta (c)} \quad\text{ if } 
\xi_{n+1}^2\le\xi_n^2.
$$
Testing $\Cal C_R^{\xi,\tau}$ on $g=\chi_{[-R,R]}$ implies that
$\|\Cal C_R^{\xi,\tau}\|\ge c \log R$ if 
$|\xi_{n+1}|\le |\xi_n|$ and the asserted lower bound on $\Cal C_R$ is proved.
Thus the better bound of 
Corollary 3.1 indicates a strictly
noncommutative phenomenon.

\endremark

\subhead{\bf A first decomposition}
\endsubhead
In view of the product type singularity of the kernel it is
natural to introduce a dyadic decomposition in the variables
$x_2-y_2$ and $x_1-y_1$ (if the latter is large). For this let
$\eta_0$ be a smooth nonnegative even function on the real line so
that $\eta_0(s)=1$ if $|s|\le 1/2$ and $\eta_0(s)=0$ if $|s|\ge
3/4$. We also assume that $\eta_0'$ has only a finite number of
sign changes. Let
$$
\eta(s)=\eta_0(s/2)-\eta_0(s)
$$
so that $\eta$ is supported in $[1/2,3/2]\cup[-3/2,-1/2]$. For
pairs of integers $j=(j_1,j_2)\in \Z^2$, with $j_1> 0$ let
$$
\chi_j(x)=\chi_{1,j_1}(x_1)\chi_{2,j_2}(x_2)=\frac{2^{j_1}}{\beta(x_1)}
\eta(2^{-j_1}x_1)
\frac{2^{j_2}}{x_2}
\eta( 2^{-j_2}x_2).
\tag 3.4
$$
In particular $\chi_j$ has the cancellation property
$$
\int\chi_j(x_1,x_2)dx_2=0 \quad\text{ for all } x_1.
$$
It will sometimes be useful (see  Section~4 below) to use the
cut-off function
$$
\widetilde \chi_j(x)=\chi_j(x)\sgn(x_2),
\tag 3.5
$$
together with the relation
$$
\chi_j(x-y)\sin( \theta(x,y))= \widetilde\chi_j(x-y)\sin
(|\theta(x,y)|),
\tag 3.6
$$
which follows from the evenness of the function $t\mapsto
t^{-1}\sin(At)$ and the positivity of $\beta$.

Let
$$
T_j(x,y)=2^{-j_1-j_2} \chi_j(x-y)
e^{i\gamma\Psi(x,y)}\sin\theta(x,y);
\tag 3.7
$$
then we wish to estimate the $L^2$ operator norm of
$$
\cT^R f(x)=\sum\Sb 10< j_1\le \log R\\j_2\in\Z\endSb \int
T_{j}(x,y)f(y) dy.
\tag 3.8
$$

\subhead{\bf Preliminary estimates}
\endsubhead
We shall now verify that the operator norm of $\cO^R- \cT^R$ is
uniformly bounded.  To this end, we consider, for fixed
$(x_1,y_1)$, the operator $\cB^{x_1,y_1}$ acting on functions in
$C^\infty_0(\R)$, which has the distribution kernel
$$
B^{x_1,y_1}(x_2,y_2)= e^{i\gamma\Psi(x,y)}
\sin\theta(x,y) (x_2-y_2)^{-1}.
\tag 3.9
$$

\proclaim{Lemma 3.2}
For each $(x_1,y_1)$ the operator $\cB^{x_1,y_1}$ extends to a
bounded operator on $L^2(\R)$ with norm bounded independently of
$(x_1,y_1)$.
\endproclaim
\demo{\bf Proof}
For $\ep,\ep'\in\{\pm 1\}$, define
$$
E^\ep_{x_1y_1}=\{(s,t): \ep(s+t+p'(x_1+y_1))\ge 0\}.
$$
One computes that
$$
2i \cB^{x_1,y_1}g(s) =\sum_{\ep'=\pm 1}\sum_{\ep=\pm 1} \ep'
e^{i\rho_\ep(s)} \int\chi_{E^\ep_{x_1,y_1}}(s,t) g(t)
e^{i\sigma_\ep(t)} \frac {dt}{s-t},
$$
where
$$
\align
\rho_\ep(s)&\equiv \rho_{\ep,x_1,y_1}(s)= (x_1-y_1) s^2-sp(x_1+y_1)+\ep
\beta(x_1-y_1)(s^2+p'(x_1+y_1)s)
\\
\sigma_\ep(t)&\equiv \sigma_{\ep,x_1,y_1}(t)= (x_1-y_1) t^2+tp(x_1+y_1)-\ep
\beta(x_1-y_1)(t^2+p'(x_1+y_1)t).
\endalign
$$
The uniform boundedness of $\cB^{x_1y_1}$ on $L^2(\R)$ follows
from the boundedness of Hilbert transforms and Hilbert
integrals.\qed
\enddemo

\proclaim{Lemma 3.3}
Let $\cE$ be an operator bounded on $L^2(\R)$, with nonnegative
kernel $k(s,t)$. Let
$$
\cS f(x)=\int k(x_1,y_1) |\cB^{x_1,y_1} [f(y_1,\cdot)](x_2)|dy_1.
$$
Then $\cS$ is bounded on $L^2(\R^2)$ with operator norm $\lc
\|\cE\|_{L^2(\R)\to L^2(\R)}$.
\endproclaim
\demo{\bf Proof}
By Lemma 3.2 we have $\|\cB^{x_1,y_1}\|\le C_0$ and therefore
$$
\align
\|\cS f\|&\le
\Bigl(\int \Bigl(\int k(x_1,y_1)\|
\cB^{x_1y_1}[f(y_1,\cdot)]\| dy_1\Bigr)^2dx_1\Bigr)^{1/2}
\\
&\le C_0
\Bigl(\int \Bigl(\int k(x_1,y_1) \bigl(\int|f(y_1,x_2)|^2
dx_2\bigr)^{1/2} dy_1\Bigr)^2 dx_1\Bigr)^{1/2},
\endalign
$$
and the result follows from the assumed $L^2$ boundedness of the
operator $\cE$ acting on the function $y_1\mapsto
\|f(y_1,\cdot)\|_{L^2(\R)}$. \qed
\enddemo

\proclaim{Lemma 3.4}
The operator $\cO^R-\cT^R$ is bounded on $L^2(\R^2)$ with an
admissible operator norm uniformly in $R$.
\endproclaim

\demo{\bf Proof}
Let $E_R=[-40,40]\cup[R/4,4R]\cup [-4R,-R/4]$ and
$k_R(s)=\chi_{E_R}(s)\,\beta(s)^{-1}$. Observe that the $L^1(\R)$
norm of $k_R$ is uniformly bounded in $R$. Note that
$$
|\cO^R f(x)- \cT^Rf(x)|\lc
\int k_R(x_1-y_1)\,\big|\cB^{x_1y_1}[f(y_1,\cdot)](x_2)\big| dy_1
$$
so that the assertion follows from Lemma 3.3.\qed
\enddemo

By Lemma 3.4 it suffices to show the bound
$$
\|\cT^R\| =O(\log\log R)
\tag 3.10
$$
for large $R$. The next four sections will be devoted to the proof
of (3.10). The argument relies on a crucial cancellation property
for the affine case, where $p(x)=ax+b$, for which one obtains the
bound $\|\cT^R\|=O(1)$. This will be carried out in Section~4.
 The general case involves an approximation by operators
which share the properties of the affine case; for various
remainder terms one uses the oscillatory properties of the phase
function and Hilbert integral arguments. The basic decomposition
describing the remainder terms and relevant orthogonality
arguments is introduced in Section~5; here we state several
propositions containing estimates for the constituents in the
basic decomposition and deduce the main estimate (3.10). Section~6
contains a few auxiliary facts and Section~7 contains the proof of
the propositions.

\head{\bf 4. Boundedness for affine polynomials}
\endhead

Let $I$ be a set of pairs $(j_1,j_2)$ with the property that
$j_1,j_2\in \Z$ and $j_1\ge 10$. Define
$$
\cT f(x)= \sum_{j=(j_1,j_2)\in I}
2^{-j_1-j_2}
\int  e^{i\gamma\Psi(x,y)}\sin\theta(x,y)\chi_j(x-y) f(y)\,dy
\tag 4.1
$$

\proclaim{Theorem 4.1}
Assume that $\alpha_0, \alpha_1 \in\R$, and that
$$
p(s) =\alpha_1 s+\alpha_0.
$$
Suppose that $1/2\le |\gamma|\le \Gamma$. Then the operator $\cT$
extends to a bounded operator on $L^2(\R^2)$, and
$$
\|\cT\|\le C_{\Gamma},
$$
where $C_\Gamma$ does not depend on $I$, $\alpha_0$, or
$\alpha_1$.
\endproclaim

\demo{\bf Proof}
We have now
$$
\aligned
\Psi(x,y)&=(x_1-y_1)(x_2^2+y_2^2)-(x_2-y_2)(\alpha_1(x_1+y_1)+\alpha_0)
\\
\theta(x,y)&= \beta(x_1-y_1)|x_2+y_2+\alpha_1|(x_2-y_2)
\endaligned
\tag 4.2
$$
and, setting
$A(x)=\tfrac{\alpha_1^2}2 x_1-2\alpha_1x_1x_2-\alpha_0 x_2$,
we compute that
$$
\Psi(x_1,x_2-\tfrac{\alpha_1}2,y_1,y_2-\tfrac{\alpha_1}2)
=(x_1-y_1)(x_2^2+y_2^2)
+A(x)-A(y);
$$
moreover
$$
\theta(x_1,x_2-\tfrac{\alpha_1}2,y_1,y_2-\tfrac{\alpha_1}2)=
\beta(x_1-y_1)|x_2+y_2|(x_2-y_2).
$$
{}From (4.1) and (4.2), we see that we can reduce matters to the
case where $p=0$, after a translation in the $x_2$ variable and a
conjugation with a multiplication operator of norm $1$. Therefore
we shall now work with (4.2) where $\alpha_0=\alpha_1=0$, and
consider the integral operator $\cK$ with kernel
$$
K= \sum_{j\in I} K_j,
$$
where
$$
K_j(x,y)= 2^{-j_1-j_2} e^{i\gamma(x_2^2+y_2^2)(x_1-y_1)}
\sin(\beta(x_1-y_1)|x_2^2-y_2^2|)\widetilde \chi_j(x-y);
$$
see  formula (3.6). For $\xi_1\in \R$, let
$$
S_{\xi_1}g(x_2)
=\sum_{j_2} \int g(y_2) h(\xi_1,x_2,y_2,j_2)  2^{-j_2} dy_2,
$$
where
$$
h(\xi_1,x_2,y_2,j_2)=\widetilde\chi_{2,j_2}(x_2-y_2)
\sum_{j_1:(j_1,j_2)\in I} h_{j_1}(\xi_1,x_2,y_2),
$$
with
$$
h_{j_1}(\xi_1,x_2,y_2)
= 2^{-j_1}
\int e^{i\sigma\gamma(x_2^2+y_2^2)-i\sigma\xi_1}
\sin(\beta(\sigma)|x_2^2-y_2^2|)
\chi_{1,j_1}(\sigma) d\sigma.
$$
Then
$$
\widetilde{\cT f}(\xi_1,x_2)= S_{\xi_1}[\widetilde
f(\xi_1,\cdot)](x_2),
$$
where $\widetilde{f }$ denotes the Fourier transform of $f$ with
respect to the first variable. Thus it suffices to fix $\xi_1$ and
show that $S_{\xi_1}$ is bounded on $L^2(\R)$ uniformly in
$\xi_1$.

\proclaim{Lemma 4.2}

(i) There is a constant $C$ so that
$$
|h(\xi_1,x_2,y_2,j_2)|\le C
\tag 4.3
$$
for all $\xi_1$, $x_2$, $y_2$, $j_2$.  Moreover
$$
h(\xi_1,x_2,y_2,j_2)=0 \quad \text{ if }
|x_2-y_2|\notin[2^{j_2-2}, 2^{j_2+2}] .
\tag 4.4
$$

(ii)  For each $j_1$
$$
|h_{j_1}(\xi_1,x_2,y_2)|\lc 2^{j_1}|x_2^2-y_2^2|.
\tag 4.5
$$

(iii) Suppose that $|\xi_1-\gamma(x_2^2+y_2^2)|\ge 2
|x_2^2-y_2^2|$. Then
$$
|h_{j_1}(\xi_1,x_2,y_2)|\le C_N
(2^{j_1}|\xi_1-\gamma(x_2^2+ y_2^2)|)^{-N}.
\tag 4.6
$$
\endproclaim

\demo{Proof} The assertion (ii) follows immediately from the
inequality $|\sin \alpha|\le |\alpha|$. Moreover (4.4) is
immediate from the definitions.
In what follows we shall use  simple properties of $\beta$ stated in (6.1), (6.2) below.

We now prove the uniform boundedness of $h$. Since $\chi_{1,j_1}$
is an even function,
$$
\align
h_{j_1}(\xi_1,x_2,y_2)&=
2^{-j_1 +1}\int_{\sigma>0} \cos(\sigma(\xi_1-\gamma (x_2^2+ y_2^2)))
\sin (\beta(\sigma)|x_2^2-y_2^2|)\chi_{1,j_1}(\sigma)d\sigma
\\
&=h_{j_1}^+(\xi_1,x_2,y_2)- h_{j_1}^-(\xi_1,x_2,y_2),
\endalign
$$
where
$$
h_{j_1}^\pm(\xi_1,x_2,y_2)= 2^{-j_1}\int_{\sigma>0}
\sin(\phi^{\pm}(\sigma;\xi_1,x_2,y_2))\chi_{1,j_1}(\sigma) d\sigma
,
$$
with
$$
\phi^{\pm}(\sigma;\xi_1,x_2,y_2)=
\sigma(\xi_1-\gamma (x_2^2+ y_2^2))
\pm\beta(\sigma)|x_2^2-y_2^2|.
$$

Observe that
$$
\align
(\phi^\pm)'(\sigma)&=\xi_1-\gamma (x_2^2+ y_2^2)
\pm\beta'(\sigma)|x_2^2-y_2^2|
\\
(\phi^\pm)''(\sigma)&=\pm\beta''(\sigma)|x_2^2-y_2^2| ,
\endalign
$$
so that
$$
|(\phi^\pm)''(\sigma)|\approx 2^{-3j_1}|x_2^2-y_2^2|
$$
in the support of $\chi_{1,j_1}$ 
(see
  (6.1) below). Using van der
Corput's Lemma (\cite{33, p.~334}), we obtain the inequality
$$
|h_{j_1}^\pm(\xi_1,x_2,y_2)|\le C 2^{j_1/2}|x_2^2-y^2_2|^{-1/2}.
\tag 4.7
$$
Now since $\sigma>0$, we have in view of (6.2.1) below
$$
\sin\phi^\pm(\sigma)= \sin\bigl(\sigma(\xi_1-\gamma(x_2^2+y_2^2) \pm\frac
12|x_2^2-y_2^2|)\bigr)+ O(|x_2^2-y_2^2|/ \sigma),
$$
so that
$$
\aligned
h_{j_1}^\pm(\xi_1,x_2,y_2)&= 2^{-j_1}\int_{\sigma>0}
\sin\bigl(\sigma(\xi_1-\gamma(x_2^2+y_2^2)\pm\frac 12
|x_2^2-y_2^2|)\bigr) \chi_{1,j_1}(\sigma) d\sigma \\&\qquad+
r_{j_1}^\pm(\xi_1,x_2,y_2),
\endaligned
\tag 4.8
$$
where the error terms $r_{j_1}^\pm$ satisfy the estimate
$$
|r_{j_1}^\pm(\xi_1,x_2,y_2)|\lc 2^{-j_1}|x_2^2-y_2^2|.
\tag 4.9
$$
Concerning the integral in (4.8),
observe that
$$
\big| \sum_{j_1\in \Cal E} \int_{\sigma>0} \sin(A\sigma) 2^{-j_1}\chi_{1,j_1}(\sigma) d\sigma
\big|\le C,
\tag 4.10
$$
where the sum is over any finite set $\Cal E$ consisting  of positive $j_1$; the
constant $C$ can be chosen independently of $\Cal E$ and of $A$.
To see this we use the inequality $|\sin \alpha|\le |\alpha|$ for
the terms with $A2^{j_1}\le 1$ and integration by parts for the
terms with $A2^{j_1}> 1$. {}From (4.10),
$$
\big|\sum_{j_1\in \Cal E} \bigl(h_{j_1}^\pm(\xi_1,x_2,y_2)-
r_{j_1}^\pm(\xi_1,x_2,y_2)\bigr)\big|\le C,
\tag 4.11
$$
where $\Cal E$ is again any set of positive indices  and the bound is
uniform in $\xi_1, x_2, y_2, j_2$.

Now an application of formulae (4.7), (4.8), (4.9) and (4.11)
shows that $h$ is uniformly bounded. Finally, the estimate (iii)
follows by integration by parts, using the lower bound
$$
|(\phi^\pm)'(\sigma)|\ge \frac 12|\xi_1-\gamma(x_2^2+y_2^2)|,
$$
in the present case of assertion (iii), 
see
 formula (6.2.2) below.
Moreover, if $\nu\ge 2$, then
$$
|\partial_\sigma^\nu\phi^\pm(\sigma)|=
|\beta^{(\nu)}(\sigma)||x_2^2-y_2^2| \lc
(1+|\sigma|)^{-|\nu|-1}|x_2^2-y_2^2|,
$$
which is an acceptable upper bound in (iii). \qed\enddemo

In what follows, $\xi_1$ will be fixed, and we shall not always
indicate the dependence of the operators on $\xi_1$. For $M\in \Z$,
let
$$
\cS^M_{j_2} g(x_2)=
\eta(2^{-M}x_2) 2^{-j_2}\int g(y_2)
h(\xi_1,x_2,y_2,j_2) dy_2.
$$

Let $C_0$ be an integer with
$2^{C_0-100}\ge \Gamma$.
We split
$$
\cS_{\xi_1} = \cS+\sum\Sb (j_2,M)\\M\le j_2+C_0\endSb \cS^M_{j_2} .
$$
It is easy to see using the uniform boundedness of $h$ and the
definition of the cut-off functions that
$$
\sum\Sb (j_2,M)\\M\le j_2+C_0\endSb
| \cS^M_{j_2}g(x_2)|
\lc \int_{|x_2|\le 2^{C_0+2}|x_2-y_2|}\frac{|g(y_2)|}{|x_2-y_2|} dy_2.
\tag 4.12
$$
The integral on the right hand side in (4.12) is a standard
Hilbert integral and therefore defines a  bounded operator on $L^2(\R)$ (see \cite{32, p.~271}).

We let
$$
\cS^M=\sum_{ j_2< M-C_0} \cS^M_{j_2}
$$
and the kernel $S^M(x_2,y_2)$ is supported where $2^{M-1}\le
|x_2|\le 2^{M+1}$, $2^{M-2}\le |y_2|\le 2^{M+2}$. Therefore the
almost orthogonality property
$$
\|\cS\|\lc \sup_M\|\cS^M\|
\tag 4.13
$$
holds. Thus it suffices to prove a uniform estimate for the
operators $\cS^M$. We split $\cS^M=\cP^M+ (\cS^M-\cP^M)$, where
$$
\cP^M g(x_2)= \sum_{j_2<M-C_0}\eta_0( 2^{-M-j_2-C_0+10}(\xi_1-2\gamma x_2^2))
\cS^M_{j_2}g(x_2).
$$
We first show that the operators $\cP^M$ are uniformly bounded.
Since $j_2\le M-C_0$, we observe that the conditions
$|2^{-M-j_2-C_0 +10}(\xi_1-2\gamma x_2^2))|\le 1$ and
$2^{M-1}\le|x_2|\le 2^{M+1}$ imply that $\xi_1(2\gamma)^{-1} >0$
and $\xi_1(2\gamma)^{-1} \approx 2^{2M}$, and therefore
$$
\big||x_2|-(\xi_1/2\gamma)^{1/2}\big|\le C_1|x_2-y_2|.
$$
Consequently
$$
|\cP^M g(x_2)|
\lc \int_{||x_2|-(|\xi_1/2\gamma|)^{1/2}|
\le C_1|x_2-y_2|}
\frac{|g(y_2)|}{|x_2-y_2|} dy_2.
$$
The right hand side is a sum of two operators, each of them a
Hilbert integral operator composed with translation operators.
Therefore it defines a bounded operator on $L^2(\R)$ and
$$
\|\cP^M\|=O(1).
$$

Next we consider the operator $\cS^M-\cP^M$ which we split as
$$
\align
[\cS^M-\cP^M]g(x_2)
&=\sum_{j_2<M-C_0}(1-\eta_0(2^{-M-j_2-C_0+10}(\xi_1-2\gamma x_2^2)))
\cS^M_{j_2} g(x_2)
\\
&=\sum_{r\in\Z}\eta(2^{-r}(\xi_1-2\gamma x_2^2))
 \sum\Sb j_2<r-M-C_0+10\\ j_2<M-C_0\endSb
\cS^M_{j_2}g(x_2)
\\&=\sum_{r\in \Z} \cQ^M_rg(x_2),
\endalign
$$
say.  Since $2^{C_0}\ge 2^{100}\Gamma$, we have
$|2\gamma(x_2^2-y_2^2)|\le 2^{M+j_2+3}|\Gamma|\le
2^{r-C_0+13}|\Gamma| \le 2^{r-10}$, and hence
$$
|\xi_1-2\gamma y_2^2|=
|\xi_1-2\gamma x_2^2+2\gamma(x_2^2-y_2^2)|\approx 2^r.
\tag 4.14
$$
Thus $|\xi_1-2\gamma x_2^2|\in( 2^{r-1}, 2^{r+1})$, which implies
that $|\xi_1-2\gamma y_2^2|\approx 2^{r}$ and we can deduce the
almost orthogonality property
$$
\|\cS^M-\cP^M\|\lc\sup_r\|\cQ^M_r\|.
\tag 4.15
$$

Now, analogously to (4.14), we also have
$$
|\xi_1-\gamma (x_2^2+y_2^2)|\approx 2^r.
\tag 4.16
$$
By Lemma 4.2 (ii) and (iii),
$$
|h_{j_1} (\xi_1,x_2,y_2)|\lc \min \{ 2^{-j_1-r},
2^{j_1+j_2+M}\}
$$
if $|x_2-y_2|\approx 2^{j_2}$ and $|x_2|\approx |y_2| \approx
2^M$, and $|\xi_1-\gamma(x_2^2+y_2^2)|\approx 2^r\gg 2^{M+j_2}$.
Therefore
$$
|h(\xi_1,x_2,y_2,j_2)|\lc 2^{(M+j_2-r)/2},
$$
and it follows that
$$
\|\cQ^M_r\|\lc \sum_{M+j_2\le r}
2^{(M+j_2-r)/2}\le C.
$$
This now implies the uniform boundedness of the operators
$\cS^M-\cP^M$. Together with the $L^2$ boundedness of $\cP^M$ and
the orthogonality property of the operators $\cS^M$ this completes
the proof of Theorem 4.1. \qed
\enddemo

\head{\bf 
5. Basic decompositions and outline of the proof for $\boldkey
n\boldsymbol \ge \boldkey 2$}
\endhead

We shall now assume that $n\ge 2$ and that $p$ is a nonaffine
polynomial of degree $\le n$. Since we are estimating the operator
$\cT^R$ we shall assume that sums in $j$ are always taken over
subsets of $\{(j_1,j_2):10<j_1\le \log R\}$.

We begin by refining the dyadic decomposition from Section~3.
Using the cut-off functions $\eta_0$ and $\eta$ as defined in
Section~3 we set

$$
a_m(\sigma)= \prod_{\nu=2}^{\deg(p)-1}\eta_0\bigl( 2^{m+10} \frac
{p^{(\nu+1)}(\sigma)}{p^{(\nu)}(\sigma)}\bigr)
\tag 5.1
$$
if $\deg(p)\ge 3$ and
$a_m(\sigma)=1$ if $\deg(p)=2$. Next,
$$
b_m(X_1, X_2)=\eta_0\bigl(2^{m+10}
\frac {p''(X_1)}{p'(X_1)+X_2}\bigr).
\tag 5.2
$$
Moreover, let
$$
\align
h_{l}(X_1,X_2)&= \eta_0(2^{-l-10}(X_2+p'(X_1)),
\tag 5.3.1
\\
h_{l,r}(X_1,X_2)&= \eta(2^{-l+r-10}(X_2+p'(X_1)),
\tag 5.3.2
\endalign
$$
so that $h_l=\sum_{r>0}h_{l,r}$ a.e.

Now let $T_j(x,y)$ be as in (3.7); our basic splitting (assuming $j_1>10$) is
$$
T_j=H_j+U_j+W_j+\sum_{r>0}V_j^r,
\tag 5.4
$$
where
$$
\align
H_j(x,y)&=T_j(x,y)(1-a_{j_1}(x_1+y_1)),
\tag 5.5.1
\\
U_j(x,y)&=T_j(x,y)a_{j_1}(x_1+y_1)(1-b_{j_1}(x+y)),
\tag 5.5.2
\\
V_j^r(x,y)&=T_j(x,y)a_{j_1}(x_1+y_1)b_{j_1}(x+y) h_{j_2,r}(x+y),
\tag 5.5.3
\\
W_j(x,y)&=T_j(x,y)a_{j_1}(x_1+y_1)b_{j_1}(x+y)(1- h_{j_2}(x+y)).
\tag 5.5.4
\endalign
$$
Let $\cH_j$, $\cU_j$, $\cV_j^r$, $\cW_j$ be the corresponding
operators. Let $\cH$, $\cU$, $\cV^r$, $\cW$ denote the operators
$\sum_j \cH_j$, $\sum_j \cU_j$, $\sum_j \cV_j^r$, and $\sum_j
\cW_j$.

We shall also use the notation
$$
\align
u_j(x,y)&=a_{j_1}(x_1+y_1)(1-b_{j_1}(x+y)),
\tag 5.6.1
\\
v_{j}^r(x,y)&=a_{j_1}(x_1+y_1)b_{j_1}(x+y) h_{j_2,r}(x+y),
\tag 5.6.2
\\
w_j(x,y)&=a_{j_1}(x_1+y_1)b_{j_1}(x+y)(1- h_{j_2}(x+y)).
\tag 5.6.3
\endalign
$$

\proclaim{Proposition 5.1 }
The operator $\cH$ is bounded on $L^2(\R^2)$.
\endproclaim

\proclaim{Proposition 5.2}
Let $\cU_j^L$ be the operator with kernel~$U_j^L$ given by
$$
U_j^L(x,y)= U_j(x,y)\eta(2^{-L}p''(2x_1)),
$$
and let $\cU^L=\sum \cU_j^L$. Then

(i)
$$\|\cU\|\lc \sup_L\|\cU^L\|.
\tag 5.7
$$

(ii)
$$ \|\cU_j^L\|\lc \min\{ 2^{L+2j_1+j_2}, 2^{-(L+2j_1+j_2)/4}\}.
\tag 5.8
$$

(iii)
$$\|(\cU^L_j)^*\cU^L_k\|
+\|\cU^L_j(\cU^L_k)^*\|\lc 2^{-|j_1-k_1|/2}.
\tag 5.9
$$
\endproclaim

\proclaim{Proposition 5.3}

(i)
$$
\|(\cV_j^r)^*\cV_k^r\|
+\|\cV_j^r(\cV_k^r)^*\|\lc 2^{-r-|j_2-k_2|/2}.
\tag 5.10
$$

(ii)
$$
 \|\cV_j^r\|\lc \min\{ 2^{2j_2+j_1 - r},   2^{r/2-(2j_2+j_1)/4}\}.
\tag 5.11
$$
\endproclaim

\proclaim{Proposition 5.4}
For $M\in \Z$, $L\in \Z$, let $\cW^{M,L}_j$ be the operator with
kernel
$$
W^{M,L}_j(x,y)=W_j(x,y)\eta(2^{-M}(2x_2+p'(2x_1)))
\eta(2^{-L}p''(2x_1)),
\tag 5.12
$$
and let $\cW^{M,L}= \sum_j \cW^{M,L}_j$. Then

(i)
$$
\|\cW\|\lc \sup_{M,L}\|\cW^{M,L}\|.
\tag 5.13
$$

(ii)
$$
 \|\cW_j^{M,L}\|\lc \min\{ 2^{M+j_1+j_2},
2^{-(M+j_1+j_2)/4}\}.
\tag 5.14
$$

(iii) The estimate (5.14) also holds if $W_j^{M,L}(x,y)$ is
replaced by $W_j^{M,L}(x,y) \rho_j(x,y)$ where $\rho_j$ satisfies
$\partial_x^\alpha\rho_j,\,
\partial_y^{\alpha} \rho_j =
O(2^{-j_1|\alpha_1|-j_2|\alpha_2|})$, for $\alpha_1,\alpha_2\in \{0,1\}$.
\endproclaim

The previous propositions are enough to obtain a uniform bound on
the operators $\cU$ and $\cV$. For $\cW$, an analogue of the
crucial orthogonality properties (5.9) and (5.10) is missing, and
we shall instead use an approximation by operators treated in
Section~4.

\proclaim{Proposition 5.5}
Suppose that $m\ge 0$. Fix $M\in \Z$ and $L\in \Z$, and let
$I$ be a set of integer pairs $j=(j_1,j_2)$ satisfying
$$
\gathered
M-L-m(1+\frac 1{2n})<j_1\le M-L-m
\\
L+2j_1+j_2\le 0
\\10\le j_1<\log R.
\endgathered
\tag 5.15
$$
Let
$$
\cZ^{I}
 =\sum_{j\in I} \cW^{M,L}_j.
$$
Then $\cZ^{I}$ is bounded on $L^2(\R^2)$ and
$$
\|\cZ^{I}\|_{L^2\to L^2} \le C,
\tag 5.16
$$
where the admissible constant $C$ is independent of $I,L,M,m,R$.
\endproclaim

Taking Propositions 5.1--5.5 for granted, we are now able to give a
proof of the main theorem.

\demo{\bf Proof of Theorem 3.1}
By the discussion in Section~3 it suffices to prove the estimate
(3.10). In view of Proposition~5.1, we have to bound $\cU$, $\cV$
and $\cW$. In order to bound $\cU$, it is sufficient to obtain a
uniform bound for the operators $\cU^L$, by (5.7). Let
$$
\cU^{L,\ell}=\sum \Sb L+2j_1+j_2=\ell\endSb
\cU^L_j.
$$
Suppose $L+2j_1+j_2=\ell$, $L+2k_1+k_2=\ell$ and $|j_1-k_1|=s$; then
by (5.8) and (5.9),
$$
\|(\cU^L_j)^*\cU^L_k\|
+\|\cU^L_j(\cU^L_k)^*\|\lc \min\{2^{2\ell}, 2^{-\ell/2}, 2^{-s/2}\},
$$
and by the Cotlar--Stein Lemma (\cite{33, p.~280})
it follows that
$$
\cU^{L,\ell}\lc \sum_{s=0}^\infty
\min\{2^{\ell}, 2^{-\ell/4}, 2^{-s/4}\}
\lc (1+|\ell|)\min\{2^\ell, 2^{-\ell/4}\}.
$$
Summing over $\ell$ yields the desired uniform bound for $\cU^L$
and thus the boundedness of $\cU$.

The operator $\cV^r$ is handled similarly.
Now let
$\cV^{\ell,r}=\sum_{2j_2+j_1=\ell}\cV_j^r$. {}From Proposition 5.3, we have
$$
\|(\cV_j^r)^*\cV_k^r\|+
\|(\cV_j^r(\cV_k^r)^*\|\lc \min\{2^{2\ell - 2r}, 2^{r-\ell/2}, 2^{-r-s/2}\}
$$
if $2j_2+j_1=2k_2+k_1 = \ell, |j_2-k_2|=s$,
and we obtain from the Cotlar--Stein Lemma that
$$
\|\cV^{\ell,r}\|\lc
\sum_{s\ge 0}\min\{2^{\ell - r}, 2^{r/2-\ell/4}, 2^{-r/2-s/4}\}
$$
and thus
$$
\align
\|\cV\|\le& \sum_{r>0}\sum_{\ell \in \bbZ}\|\cV^{\ell,r}\|
\\ \lc&
\sum_{r> 0} \sum_{\ell\le 0} 2^{-|\ell|-r}(1+|\ell|+r)
+
\sum_{r>0} \sum_{0\le \ell\le 4r} 2^{-r/2}
\\& +
\sum_{r> 0} \sum_{ \ell> 4r} \Bigl( 2^{-r/2-(\ell-4r)/4}+ 2^{-r/2}2^{r-\ell/4}(\ell-4r)\Bigr),
\endalign
$$
and  $\cV$ is easily seen to be bounded on $L^2$.

Now we turn to the operator $\cW$.  By (5.13) it suffices to obtain
a uniform bound for $\cW^{M,L}$.
We note that $\cW_j^{M,L}=0$ if $L+j_1\ge M$.
Therefore by (5.14)
$$
\sum\Sb L+2j_1+j_2\ge 0\endSb\bigl\|\cW^{M,L}_j\bigr\|
\lc
\sum\Sb L+2j_1+j_2\ge 0\\L+j_1\le M\endSb
2^{-(M+j_1+j_2)/4}\le C.
$$
For sums of terms $\cW^{M,L}_j$ which satisfy $L+2j_1+j_2< 0$ we
use Proposition 5.5. For $s=1,2,\dots$, let
$$
I^{M,L}_{s,R}=\{j:M-L-(\tfrac{2n+1}{2n})^s<j_1\le M-L-
(\tfrac{2n+1}{2n})^{s-1}, L+2j_1+j_2<0, 10<j_1\le \log (10+R)\}.
$$
Then from Proposition 5.5,
$$
\Bigl\|
\sum\Sb j\in I^{M,L}_{s,R}\endSb
\cW^{M,L}_j\Bigr\| \le C,
$$
uniformly in $s$, $R$ and $M$, $L$. Now for fixed $M,L,R$ the sets
$I^{M,L}_{s,R}$ are nonempty for no more than $C_0\log\log(10+ R)$
choices of $s$; here $C_0$ is admissible. Summing over $s$ we
 see that $\|\cW^{M,L}\|=O(\log\log(10+R))$, with an admissible constant,
and by (5.13) we obtain the same bound for $\cW$.
\qed\enddemo

\head
{\bf 6. Auxiliary Lemmas}
\endhead
We first collect formulae for the derivatives of $\beta$, $\Psi$ and $\theta$.
\proclaim{Lemma 6.1}
(i)
$$
\beta'(s)=\frac{s}{4\beta(s)},\qquad
\beta''(s)=\frac 1{4[\beta(s)]^3}
\tag 6.1
$$
and
$$
\gather
\Big|\beta(s)-\frac{|s|}2\Big|
=\frac{2}{2\beta(s)+|s|}\le \frac 1{|s|}
\tag 6.2.1
\\
\Big|\beta'(s)-\frac{1}2 \sgn(s)\Big|\le \frac{4}{4+s^2}.
\tag 6.2.2
\endgather
$$

(ii) Let $\Xi(x,y)=x_2+y_2+p'(x_1+y_1)$.
Suppose that $\Xi(x,y)\neq 0$. Then

$$
\align
\Psi_{x_1}(x,y) &= 2 x_2^2-(x_2-y_2)\Xi(x,y)
\tag 6.3\\
\Psi_{x_2}(x,y) &= 2 x_2(x_1-y_1)- p(x_1+y_1)
\tag 6.4
\\
\theta_{x_1}(x,y)&= (x_2-y_2)\bigl[ \beta'(x_1-y_1)|\Xi(x,y)|
+\beta(x_1-y_1)p''(x_1+y_1)\sgn(\Xi(x,y))\bigr]
\tag 6.5
\\
\theta_{x_2}(x,y)
&=\beta(x_1-y_1)(2x_2+p'(x_1+y_1))\sgn(\Xi(x,y)).
\tag 6.6
\endalign
$$

(iii)
$$
\align
\Psi_{x_1y_1}(x,y)&=-(x_2-y_2)p''(x_1+y_1)
\tag 6.7
\\
\theta_{x_1y_1}(x,y)
&=(x_2-y_2)
\bigl[\beta(x_1-y_1)
p'''(x_1+y_1)\sgn(\Xi(x,y))-\beta''(x_1-y_1)|\Xi(x,y)|\bigr].
\tag 6.8
\endalign
$$

 (iv)
$$
\gather
\Psi_{x_1y_2}(x,y)=
2y_2+p'(x_1+y_1)
\tag 6.9
\\
\Psi_{x_2 y_1}(x,y)= -(2x_2+p'(x_1+y_1)),
\tag 6.10
\endgather
$$
$$
\align
\theta_{x_1y_2}(x,y)
&=\beta'(x_1-y_1)
\bigl[(x_2-y_2)\sgn(\Xi)-|\Xi(x,y)|\bigr]
-\beta(x_1-y_1)p''(x_1+y_1)\sgn(\Xi)
\\
&=\bigl[-
\beta(x_1-y_1)p''(x_1+y_1)-\beta'(x_1-y_1)(2y_2+p'(x_1+y_1))
\bigr]\sgn(\Xi)
\tag 6.11
\\
\theta_{x_2y_1}(x,y)&=
-\beta'(x_1-y_1)\bigl[(x_2-y_2) \sgn(\Xi)+|\Xi(x,y)|\bigr]
+\beta(x_1-y_1)p''(x_1+y_1)\sgn(\Xi)
\\&=\bigl[\beta(x_1-y_1)p''(x_1+y_1)-\beta'(x_1-y_1)(2x_2+p'(x_1+y_1))\bigr]
\sgn(\Xi).
\tag 6.12
\endalign
$$

(v)
$$
\gather
\Psi_{x_2x_2 y_1}(x,y) = -2
\tag 6.13
\\
\theta_{x_2x_2y_1}(x,y) =-2\beta'(x_1-y_1)\sgn(\Xi)
\tag 6.14
\\
\Psi_{x_2y_2}(x,y) =0
\tag 6.15
\\
\theta_{x_2y_2}(x,y) =0.
\tag 6.16
\endgather
$$
\endproclaim

\demo{\bf Proof} These are straightforward computations.\qed\enddemo

We shall now examine the properties of the cut-off functions in (5.1--5.3).
For this, the following observations are essential.

\proclaim{Lemma 6.2} Let $P$ be a polynomial, let $\ell\le \deg(P)$
and let
$$
\alpha_m(\sigma)=
\prod_{\nu=\ell}^{\deg(P)}\eta_0\bigl( 2^{m+10}
\frac {P^{(\nu)}(\sigma)}{P^{(\nu-1)}(\sigma)}\bigr).
$$

(i) Suppose that $\sigma\in \supp \alpha_m$, and
 $|\sigma-\tau|\le 2^{m+7}$.
Then for $\nu=\ell,\dots , \deg(P)$
$$
|P^{(\nu)}(\tau)-
P^{(\nu)}(\sigma)|
\le \frac 15 |P^{(\nu)}(\sigma)|
\tag 6.17
$$
and
$$
|P^{(\nu)}(\tau)|\le 2^{1-(m+4)k}|P^{(\nu-k)}(\sigma)|\quad\text{ if }
\nu-k\ge \ell.
\tag 6.18
$$

(ii) For $r=1,2,3,\dots$
$$
|\alpha_m^{(r)}(\sigma)|\le C_r 2^{-rm}.
\tag 6.19
$$
\endproclaim

\demo{\bf Proof} (i) If $\sigma=\tau$ then a slightly better
estimate than (6.18) follows from the definition of $\alpha_m$,
and then for $|\sigma- \tau|\le 2^{m+7}$ the estimate (6.18)
follows once (6.17) is proved.  To see (6.17) suppose that
$\sigma\in \supp \alpha_m$, and
 $|\sigma-\tau|\le 2^{m+7}$. Then a Taylor expansion yields
$$
\align
|P^{(\nu)}(\tau)-
P^{(\nu)}(\sigma)|&\le\sum_{k=1}^{\deg(P)-\nu}\Big|\frac {P^{(\nu+k)}(\sigma)}
{k!}(\tau-\sigma)^k\Big|
\\
&\le |P^{(\nu)}(\sigma)|\sum_{k=1}^{\deg(P)-\nu}
\frac{2^{-(m+10)k} 2^{(m+7)k}}{k!}
\\
&\le
(e^{1/8}-1)
|P^{(\nu)}(\sigma)|
\endalign
$$
and $e^{1/8}-1\le 1/5$.

(ii) follows from multiple applications of the chain rule, and the
definition of the cut-off functions.\qed
\enddemo

We now set
$$
(\tu_j,\tv_j^r,\tw_j)=(u_j\chi_j, v_j^r\chi_j, w_j\chi_j)
\tag 6.20
$$
and
$$
\align
\tu_j^L(x,y)&=\tu_j(x,y)\eta(2^{-L}(p''(2x_1)))
\tag 6.21.1
\\
\tw_{j}^{M,L}(x,y)&=
\tw_j(x,y)\eta(2^{-L}(p''(2x_1)))
\eta(2^{-M}(p'(2x_1)+2x_2))
\tag 6.21.2
\endalign
$$

\proclaim{Lemma 6.3}

For $l=1,2,3,\dots$ the following holds.

(i)

$$
|\partial_{x_1}^l\tu_j(x,y)|+
|\partial_{y_1}^l \tu_j(x,y)|
\le C_l 2^{-l j_1}.
\tag 6.22
$$

(ii)
$$
|\partial_{x_1}^l\tu_j^L(x,y)|+
|\partial_{y_1}^l \tu_j^L(x,y)|
\le C_l 2^{-l j_1}, \qquad\text{for all }j_1;
\tag 6.23
$$

(iii)
$\tw_j^{M,L}(x,y)=0$ if either $M\le j_2$ or $L+j_1\ge M$.
Moreover
$$
|\partial_{x_1}^{l_1}\partial_{x_2}^{l_2}\tw_j^{M,L}(x,y)|+
|\partial_{y_1}^{l_1}\partial_{y_2}^{l_2} \tw_j^{M,L}(x,y)|
\le C_l 2^{-l_1 j_1-l_2j_2}.
\tag 6.24
$$

(iv) For all $x_1,x_2,y_1,y_2$
$$
 \int|\partial_{x_2}\tv_j^r (x,y)|dx_2 +
\int|\partial_{y_2}\tv_j^r (x,y)|dy_2 \le C.
\tag 6.25
$$
\endproclaim

\demo{Proof} These are straightforward computations using the chain rule,
Lemma 6.2 and the definition of the cut-off functions. For (6.25),
we use the fact that the sign of $\eta_0^{\prime}$ changes finitely
many times.\qed
\enddemo

The next lemma is used to estimate various operators of Hilbert
integral type. The argument is closely related to one in \cite{9}.

\proclaim{Lemma 6.4}
Let $P$ be a polynomial of degree $\le m$. Then for $\rho>0$
$$
\iint_{|s-t|\le \rho}
\big|1-\eta_0(A\frac{P'(c_1s+c_2t)}{P(c_1s+c_2t)}\rho)\big|\, ds \,dt \le C
\frac{Am^2}{|c_1+c_2|}\rho^2.
$$
\endproclaim

\demo{Proof}
 Let
$\kappa_1<\dots<\kappa_\ell$ be the real parts of the zeroes of $P$.
For $\nu=1,\dots,\ell-1$, let
 $\mu_\nu=(\ka_\nu+\ka_{\nu+1})/2$.
Let $I_1=(-\infty,\mu_1)$,
$I_\nu=(\mu_{\nu-1},\mu_{\nu})$, $2\le \nu\le \ell-1$, and
$I_\ell=(\mu_{\ell-1},\infty)$.
Then
$$
\Big|\frac{P'(\sigma)}{P(\sigma)}\Big|\le \frac m{|\sigma-\ka_\nu|},
\text{ for }\sigma\in I_\nu.
$$
Therefore the set
$$
\{(s,t):|s-t|\le \rho, A\frac{|P'(c_1s+c_2t)|}{|P(c_1s+c_2t)|}\rho\ge 1/2\}
$$
is contained in
$$
\bigcup_{\nu=0}^{\ell}\{(s,t): (c_1s+c_2 t)\in I_\nu;\;
|c_1s+c_2 t-\ka_\nu|\le 2mA\rho;\; |s-t|\le \rho\}
$$
 which is easily
seen to be of measure $O(\rho^2)$; in particular one may check the
asserted dependence on $c_1,c_2$.\qed
\enddemo

\remark{Remark} We shall use this lemma just for the regular case where
$(c_1,c_2)=(1,1)$.
\endremark

\head{\bf 7. Proofs of Propositions 5.1-5.5}
\endhead

\subhead{7.1. Proof of Proposition 5.1}\endsubhead

We may assume that $p$ is a polynomial of degree at least
three, since otherwise $\cH=0$. For $2\le \nu\le n-1$, let
$$
Q^\nu_m(s)
=\{t:|s-t|\le 2^{m+1}, 2^{m+12} |p^{(\nu+1)}(s+t)|\ge |p^{(\nu)}(s+t)|\},
$$
and for $g\in L^2(\R)$,
$$
\cE^\nu_m g(s)= 2^{-m}\int_{Q^\nu_m(s)} g(t) \,dt
$$
and $\cE^\nu=\sum_m\cE^\nu_m$. One can use an argument in \cite{9}
to show that $\cE^\nu$ is bounded on
$L^2(\R)$. Alternatively, we use an almost orthogonality argument
based on Lemma 6.4. Specifically, denote by $k_{lm}(w,z)$
the kernel of $(\cE^\nu_l)^*\cE^\nu_m$. Since
$((\cE^\nu_l)^*\cE^\nu_m)^*= (\cE^\nu_m)^* \cE^\nu_l$ we may
assume that $l\le m$. Then, for fixed $z$,
$$
\int|k_{lm}(w,z)|dw\lc 2^{-m-l}\big|
\{(w,s):|s-w|\le 2^{l+1}, \frac{|p^{(\nu+1)}(s+w)|2^{l+1}}
{|p^{(\nu)}(s+w)|}\ge \frac 12\}\big|\lc 2^{l-m},
$$
by Lemma 6.4, and since also $\int|k_{lm}(w,z)| dz=O(1)$ for all
$w$ we see from Schur's test that
$\|(\cE^\nu_l)^*\cE^\nu_m\|=O(2^{-|m-l|/2})$. By the symmetry of
the operators $\cE^\nu_l$, one also gets
$\|\cE^\nu_l(\cE^\nu_m)^*\|=O(2^{-|m-l|/2})$ and the Cotlar--Stein
Lemma shows the $L^2$ boundedness of $\cE^\nu$. Now
$$
|\cH f(x)|\lc \sum_{\nu=2}^{n-1}\sum_{j_1} 2^{-j_1}
\int_{Q^\nu_{j_1}(x_1)}\big|\cB^{x_1,y_1}[f(y_1,\cdot)](x_2)\big| dy_1,
$$
where
$\cB^{x_1,y_1}$ is as in (3.9), and by Lemma 3.3 it follows that $\cH$ is bounded.
\qed

\subhead{\bf 7.2. Proof of Proposition 5.2}\endsubhead

Part
(i) follows from Lemma 6.2 above. Indeed suppose that
$a_{j_1}(x_1+y_1)\neq 0$,
$|x_1-y_1| \le 2^{j_1+1}$ and $2^{L-1}\le |p''(2x_1)|\le 2^{L+1}$.
Then from (6.17)
$$
|p''(2x_1)-p''(x_1+y_1)|\le \frac 15 |p''(2x_1)|\le \frac{4}{5} 2^{L-1},
$$
and similarly
$|p''(2y_1)-p''(2x_1)|\le \frac{4}{5} 2^{L-1}$. Hence if $(x,y)\in \supp U_j^L$ then 
both $p''(2x_1)$ and $p''(2y_1)$ lie in the interval $(2^{L-4}, 2^{L+2})$. This
 clearly implies that the operators $\cU^L$ are almost orthogonal
(in fact $(\cU^L)^*\cU^{L'}=0$ and
$\cU^L(\cU^{L'})^*=0$ if $|L-L'|\ge 10$).

Assuming that $L+2j_1+j_2\le 0$,
the estimate (5.8) follows from the definition of  $1-b_{j_1}$ and the
inequality
$|\sin a|\le |a|$.

Now assume that $L+2j_1+j_2\ge 1$ and write the sine as the sum of two
complex
exponentials. Then we have to estimate operators $\cR^{\ep,L}_j$
with kernels
$$
R^{\ep,L}_j(x,y)=
2^{-j_1-j_2} \chi_j(x-y) u_j^L(x,y)
e^{i{\gamma\Psi(x,y)+\ep\theta(x,y)}},
\tag 7.2.1
$$
where $\ep=\pm 1$ and $u_j^L(x,y)=u_j(x,y)\eta(2^{-L}p''(2x_1)).$

Let
$$
R_j^{\ep,L,x_2,y_2}(x_1,y_1)= 2^{j_2}R_j^{\ep,L}(x_1,x_2,y_1,y_2),
\tag 7.2.2
$$
and denote by $\cR_j^{\ep,L,x_2,y_2}$ the corresponding operator
acting on functions in $L^2(\R)$.

Let $\Phi\equiv \Phi^\ep=\gamma\Psi+\ep \theta$.
We note that
$$
|\Phi_{x_1}(x_1,y_1,x_2,y_2)-\Phi_{x_1}(x_1,z_1,x_2,y_2)|\approx
2^{L+j_2}|y_1-z_1|.
$$
This follows since by (6.7) and (6.17),
$$
|\Phi_{x_1y_1}(x_1,\tz_1,x_2,y_2)|\approx 2^{L+j_2}
$$
if $\tz_1$ is between $z_1$ and $y_1$. The derivative $
\Phi_{x_1}(x_1,y_1,x_2,y_2)-\Phi_{x_1}(x_1,z_1,x_2,y_2)$ has only
a bounded number of sign changes and we may use van der Corput's
Lemma to see that the kernel $K_{j}(y_1,z_1)$ of
$(\cR_j^{\ep,L,x_2,y_2})^* \cR_j^{\ep,L,x_2,y_2} $
 satisfies the estimate
$$
|K_{j}(y_1,z_1)|\lc 2^{-j_1}(1+2^{L+j_2+j_1}|y_1-z_1|)^{-1}.
$$
Hence it follows from Schur's test that
$$
\|(\cR_j^{\ep,L,x_2,y_2})^* \cR_j^{\ep,L,x_2,y_2} \|\lc
(L+2j_1+j_2) 2^{-L-2j_1-j_2}\lc 2^{-(L+2j_1+j_2)/2},
$$
uniformly in $x_2,y_2$. Consequently $\cR_j^{\ep,L,x_2,y_2}$ is
bounded on $L^2(\R)$, with operator norm of order at most
$2^{-(L+2j_1+j_2)/4}$, and by an averaging argument (see the proof
of Lemma 3.3) it follows that
$$
\|\cR_j^{\ep,L}\|
\lc 2^{-(L+2j_1+j_2)/4}
$$
and also that the same bound holds for $\cU^L_j$.

The orthogonality property (5.9) follows again from the argument
in Lemma 6.4.  We now give the proof for
$(\cU^L_j)^*\cU^L_k$, and without loss of generality, we may assume
that $k_1\le j_1$.

Let $K_{jk}(y,z)$ be the kernel of $(\cU^L_j)^*\cU^L_k$;
with $k_1\le j_1$.
Now for every $(x_2,z_2)$, let
$$
E^{k_1}_{x_2,z_2}=
\{(x_1,z_1): |x_1-z_1|\le 2^{k_1+1},\,
|p''(x_1+z_1)|\ge 2^{-k_1-100}|p'(x_1+z_1)+x_2+z_2|\}.
$$
By Lemma 6.4, the measure of
$E^{k_1}_{x_2,z_2}$ is $O(2^{2k_1})$. Therefore
$$
\iint|u_j^L(x,y) u_k^L(x,z)\chi_j(x-y)\chi_k(x-z)| dx_1dz_1
\le
\iint_{E^{k_1}_{x_2,z_2}} dx_1dz_1\lc 2^{2k_1}.
$$
This yields
$$
\sup_y \int |K_{jk}(y,z) | dz \lc 2^{k_1-j_1},
$$
and together with the obvious estimate
$$
\sup_z \int |K_{jk}(y,z) | dy \lc C,
$$
this implies that $\|(\cU^L_j)^*\cU^L_k\|= O(2^{(k_1-j_1)/2})$, if
$k_1\le j_1$.

For the estimation of $\|\cU^L_j(\cU^L_k)^*\|$, one uses that also
$|p''(2y_1)|\approx |p''(2x_1)|\approx 2^{L}$ on the support of
the amplitudes (as pointed out above); the argument is then the
same as for $(\cU^L_j)^*\cU^L_k$.
\qed

\subhead{\bf 7.3. Proof of Proposition 5.3}\endsubhead

We first show the  bounds asserted  for $(\cV_j^r)^*\cV_k^r$ in (5.10).
Since $\|(\cV_j^r)^*\cV_k^r\| =\|(\cV_k^r)^*\cV_j^r\|$, it suffices
to consider the case where $k_2\le j_2$. Observe that the kernel
$K_{jk}$ of $(\cV_j^r)^*\cV_k^r$ is given by
$$
\multline K_{jk}(y,z)=\\
2^{-j_1-k_1-j_2-k_2}\int
e^{i(\Psi(x,z)-\Psi(x,y))}\sin\theta(x,z)\sin\theta(x,y)
v_j^r(x,y)v_k^r(x,z)\chi_j(x-y)\chi_k(x-z) dx.
\endmultline
$$

For fixed $ x_1$ and $ z_1$, we estimate
$$
\aligned
&\iint|v_j^r(x,y) v_k^r(x,z)\chi_j(x-y)\chi_k(x-z)| dx_2dz_2
\\&\le
\iint\Sb
 |p'(x_1+z_1)+x_2+z_2|\le 2^{k_2+10-r}\\
|x_2-z_2|\le 2^{k_2+1}\endSb dx_2dz_2\lc 2^{2k_2-r},
\endaligned
\tag 7.3.1
$$
and therefore
$$
\int|K_{jk}(y,z_1,z_2)|dz_2\lc
\min\{2^{j_1}, 2^{k_1}\}2^{2k_2-r}2^{-j_1-j_2-k_1-k_2}.
$$
Now $K_{jk}$ is supported where
$|y_1-z_1|\le \max\{2^{j_1+2}, 2^{k_1+2}\}$, and so
$$
\sup_y \int |K_{jk}(y,z) | dz \lc 2^{k_2-j_2-r}.
$$
If we reverse the role of $y$ and $z$ in (7.3.1), we have to use
the less favorable bound
$$
\align
&\iint|v_j^r(x,y) v_k^r(x,z)\chi_j(x-y)\chi_k(x-z)| dx_2dy_2
\\&\le
\iint\Sb
 |p'(x_1+z_1)+x_2+z_2|\le 2^{k_2+10-r}\\
|x_2-y_2|\le 2^{j_2+1}\endSb dx_2dy_2\lc 2^{j_2+k_2-r},
\tag 7.3.2
\endalign
$$
and we obtain
$$
\sup_z \int |K_{jk}(y,z) | dy \lc C 2^{-r}.
$$
Taking the geometric mean and applying Schur's test, it follows that
$$
\|(\cV_j^r)^*\cV_k^r\|= O(2^{-r-|k_2-j_2|/2}).
\tag 7.3.3
$$
By the symmetry of $\cV_j^r$ we obtain the same bound for
$\|\cV_j^r(\cV_k^r)^*\|$.

We now turn to the assertion (ii). To obtain  the bound $\|\cV_j^r\|=O(2^{2j_2+j_1-r})$ we just
use Schur's lemma and  invoke  the estimate
 $|\sin a|\le |a|$ and the support property of $h_{j_2,r}$.

It remains to prove that $\|cV_j^r\|=O(2^{r/2-(2j_2+j_1)/4})$. Take $\ep, \ep'\in \{\pm 1\}$,
and define
$$
\Gamma_{\ep,\ep'}=\{(x,y):\sgn\gamma= \ep'\ep \,\sgn\beta'(x_1-y_1)\,
\sgn(x_2+y_2+p'(x_1+y_1))\}.
\tag 7.3.4
$$
Let $\chi_{\ep,\ep'}$ be the characteristic function of $\Gamma_{\ep,\ep'}$,
and let
$$
V_j^{r,\ep,\ep'}(x,y)
=
2^{-j_1-j_2} \chi_j(x-y)
e^{i(\gamma\Psi(x,y)+\ep\theta(x,y))}
v_{j}^r(x,y)\chi_{\ep,\ep'}(x,y),
\tag 7.3.5
$$
so that
$$
V_j^r=\sum_{\ep, \ep'\in \{-1,1\}}\frac{\ep}{2i} V_j^{r,\ep, \ep'}.
$$

It clearly suffices to prove that
$$
\|\cV_j^{r,\ep,\ep'}\|=O(2^{r/2-(2j_2+j_1)/4}).
\tag 7.3.6
$$
 The kernel $K_j(y,z)$ of
$(\cV^{r,\ep,\ep'}_j)^*\cV^{r,\ep,\ep'}_j$ is given by
$$
K_j(y,z)=
2^{-2j_1-2j_2}\int \int_{E(y,z, x_1)}
 e^{i(\Phi(x,z)-\Phi(x,y))} v_{j}^r(x,y) v_{j}^r(x,z)
\chi_j(x-y)\chi_j(x-z) \,dx_2 \,dx_1,
$$
where
$$
E(y,z,x_1)=\{x_2: (x_1,x_2,y)\in \Gamma_{\ep,\ep'},
(x_1,x_2,z)\in \Gamma_{\ep,\ep'}\}.
$$
Clearly $E(y,z,x_1)$ is the union of no more than $16$ intervals.
We note that
$$
|\Phi_{x_2x_2}(x,z)-\Phi_{x_2x_2}(x,y))|\approx|y_1-z_1|.
\tag 7.3.7
$$
To see this, apply the mean value theorem and observe that
$\Psi_{x_2x_2y_1}=-2$
and
$$
\theta_{x_2x_2y_1}=-\ep'\ep\sgn\gamma+o_{j_1},
$$
where $|o_{j_1}|\le 2^{-2j_1}$.
 Thus, since
 $|\gamma|\ge 1$ and $j_1\ge 10$ we see that
 $|\Phi_{x_2x_2y_1}|\approx 2$.
Hence we can use (7.3.7)
 to apply van der Corput's lemma
 on each of the
connected components of
$E(y,z, x_1)$. Taking into account the bound (6.25), we see that
$$
\Big|\int_{E(y,z, x_1)}
 e^{i(\Phi(x,z)-\Phi(x,y))} v_{j}^r(x,y) v_{j}^r(x,z)
 \chi_j(x-y)\chi_j(x-z) dx_2\Big|\lc 2^r |y_1-z_1|^{-1/2},
\tag 7.3.8
$$
uniformly in $x_1$, $y_2$ and $z_2$. {}From (7.3.8), it follows that
$|K_j(y,z)|$ is dominated by $2^{r-j_1-2j_2}|y_1-z_1|^{-1/2}$, and
of course it is supported where $|y_1-z_1|\le 2^{j_1+1}$,
$|y_2-z_2| \le 2^{j_2+1}$. We apply Schur's test and deduce that
$$
\|(\cV^{r,\ep,\ep'}_j)^*\cV^{r,\ep,\ep'}_j\|\lc
2^{r-(2j_2+j_1)/2};
$$
hence we get the bound (7.3.6) and consequently the bound
$\|\cV^{r}_j\|\lc
2^{r/2-(2j_2+j_1)/4}$.\qed

\subhead{\bf 7.4. Proof of Proposition 5.4}\endsubhead

Part (i) follows in view of the localization of the amplitude. Suppose
that $(x,y)\in \supp W_{j}^{M,L}$ and $\chi_j(x,y)\neq 0$. Then
$2^{M-1}\le |2x_2+p'(2x_1)|\le 2^{M+1}$ and since
$|p'(2y_1)-p'(2x_1)| \le 2^{j_1 + L + 2}$ from Lemma 6.2, we have
$$
2y_2+p'(2y_1)\in [ 2^{M-1}-2^{j_2+2} - 2^{j_1 + L+2},
2^{M+1}+2^{j_2+2} + 2^{j_1 + L+2} ];
$$
moreover the quantity
$x_2+y_2+p'(x_1+y_1)$ is also contained in this interval.  Since
$j_2\le M-10$ and $L+j_1\le M-10$, we see that
$$
\align
2^{M-2}&\le |x_2+y_2+p'(x_1+y_1)|\le 2^{M+2}
\tag 7.4.1
\\2^{M-2}&\le |2y_2+p'(2y_1)|\le 2^{M+2}.
\tag 7.4.2
\endalign
$$
Furthermore
$|x_1-y_1| \le 2^{j_1+1}$, $2^{L-1}\le |p''(2x_1)|\le 2^{L+1}$, and so
Lemma 6.2 yields $|p''(2y_1)|\approx |p''(2x_1)|\approx 2^L$. Therefore
 we can conclude that the operators
$\cW^{M,L}$ are almost orthogonal; specifically $(\cW^{M,L})^*
\cW^{M',L'} =0$ and $\cW^{M,L}(\cW^{M',L'})^*=0$ if either
$|M-M'|\ge 10$ or $|L-L'|\ge 10$. This implies (5.13). If
$M+j_1+j_2\le 0$ the estimate (5.14) follows from the fact
that $\sin a=O(|a|)$.

We now assume that $M+j_1+j_2\ge 0$.
For $\ep,\ep'\in\{\pm 1\}$, let $\chi_{\ep,\ep'}$
be the characteristic function of the set
$\Gamma_{\ep,\ep'}$, defined in (7.3.4).
Fix $L,M,\ep, \ep'$ and let
$$
\om_j(x,y)\equiv \om_j^{M,L,\ep,\ep'}(x,y):=
 \widetilde w_j^{M,L}(x,y) \chi_{\ep,\ep'}(x,y)
$$
(see (6.21.2)), and let $\cW^{M,L,\ep,\ep'}_j$ be the integral
operator with kernel
$$
W^{M,L,\ep,\ep'}_j(x,y)=
2^{-j_1-j_2} \om^{M,L,\ep,\ep'}_j(x,y)
e^{i(\gamma\Psi(x,y)+\ep\theta(x,y))}.
$$
Multiplication with the
characteristic function $\chi_{\ep,\ep'}$ does not
introduce additional singularities
in view of the localization of the symbol $w_j$; in fact,
 we have the estimates
$$
|\partial_{x}^\alpha \partial_{y}^\beta \om_j(x,y)|\le C_{\alpha,\beta}
2^{-j_1(\alpha_1+\beta_1)} 2^{-j_2(\alpha_2+\beta_2)}.
\tag 7.4.3
$$
The kernel $K_j(y,z)$ of
$(\cW^{M,L,\ep,\ep'}_j)^*\cW^{M,L,\ep,\ep'}_j$ is given by
$$
K_j(y,z)=
2^{-2j_1-2j_2}\int
 e^{i(\Phi(x,z)-\Phi(x,y))}\om_j(x,y)\om_j(x,z)
\chi_j(x-y)\chi_j(x-z) dx.
$$
In view of our assumptions that $|\gamma|\ge 1$ and $j_1\ge 10$,
we see that
$$
|\Phi_{x_2y_1}|\approx |2x_2+p'(x_1+y_1)|\approx 2^M,
$$
and also that $\Phi_{x_2y_2}\equiv 0$. Hence
$$
|\Phi_{x_2}(x,z)-\Phi_{x_2}(x,y)|\approx 2^M|y_1-z_1|.
$$
Applying van der Corput's Lemma,
$$
|K_j(y,z)|\lc
2^{-j_1-j_2}(1+2^{M+j_2}|y_1-z_1|)^{-1},
$$
 and since $K_j$ is supported where $|y_2-z_2|\le 2^{j_2+1}$,
we have
$$
\sup_y\int |K_j(y,z)| dz+\sup_z\int|K_j(y,z)|dy\lc (M+j_2+j_1)
2^{-M-j_2-j_1}.
$$
By Schur's test,
$$
\|\cW^{M,L,\ep,\ep'}_j\|\lc
(M+j_2+j_1)^{1/2} 2^{-(M+j_2+j_1)},
$$
 and this completes the proof
of (5.14).

Since we have only used property (7.4.3) our argument proves the
assertion (iii) as well.\qed

\subhead{\bf 7.5. Proof of Proposition 5.5}\endsubhead


Fix $I$ and let $\cZ\equiv \cZ^{I}$. Choose $\zeta\in C_0^\infty(\R)$,
supported in $(-1,1)$ and with the property that $\sum_{\nu\in\Z}\zeta(s-\nu)=1$.
Define the operator $\cZ_\nu\equiv \cZ^{I}_\nu$ with kernel
$$
Z_\nu(x,y)=
\zeta(2^{L+m-M+10}x_1-\nu) Z^I(x,y).
$$
In view of the localization $|x_1-y_1|\le 2^{M-L-m+1}$ we see that
the operators $\cZ_{\nu}$ are almost orthogonal; i.e.,
$(\cZ_\nu)^*\cZ_{\nu'}=0$ and $\cZ_\nu(\cZ_{\nu'})^*=0$ if
$|\nu-\nu'|\ge 100$. Therefore
$$
\|\cZ\|\lc \sup_\nu
\|\cZ_{\nu}\|.
\tag 7.5.1
$$
It hence suffices to prove a uniform estimate for the operators $\cZ_\nu$.
We wish to approximate the phase functions
$\Psi$ and $\theta$ by affine functions in the first variable.

We may suppose there is a point $c_\nu$ such that
$$
 \eta(2^{-L} p''(2c_\nu))
\neq 0 \text{ and }
|c_\nu- 2^{M-L-m}\nu|\le 2^{M-L-m-9},
 \tag 7.5.2
$$
for if not then $\cZ_\nu\equiv 0$.
Define
$$
\gather
\Psi_\nu(x,y)=(x_1-y_1)(x_2^2+y_2^2)-(x_2-y_2)p(2c_\nu)-
(x_2-y_2)p'(2c_\nu)(x_1+y_1-2c_\nu)
\\
\theta_\nu(x,y)=\beta(x_1-y_1)(x_2-y_2)|x_2+y_2+p'(2c_\nu)|.
\endgather
$$
Now $\cZ_{\nu}=\sum_{j\in A} \cZ_{\nu,j}$ where in the sum only
those $j_1$ come up which satisfy
$$
M-L-\frac{2n+1}{2n}m<j_1\le M-L-m,
$$
and the kernel of $\cZ_{\nu,j}$ is defined by
$$
\multline
Z_{\nu,j}(x,y)=2^{-j_1-j_2}
e^{i\gamma\Psi(x,y)}\sin\theta(x,y)\chi_j(x-y) \\
\times
 a_{j_1}(x_1+y_1) b_{j_1} (x+y)
\eta(2^{-L} p''(2x_1))\eta(2^{-M}(2x_2 + p'(2x_1)))
\zeta(2^{L+m-M+10}x_1-\nu) .
\endmultline
\tag 7.5.3
$$
Since we assume that $L+j_1\ll M$ the function $b_{j_1}(x+y)$ can
be omitted in (7.5.3); it is equal to $1$ on the support of the
other cut-off functions. Let
$$
z_{j}^{L,M,\nu}(x,y)=
\chi_j(x-y) \eta(2^{-L} p''(2x_1))\eta(2^{-M}(2x_2 + p'(2x_1)))
\zeta(2^{L+m-M+10}x_1-\nu) .
$$
We split $Z_{\nu,j}(x,y)$ as $\sum_{i=1}^3 Z_{\nu,j}^i(x,y)$,
where
$$
\align
&Z_{\nu,j}^1(x,y)=2^{-j_1-j_2}
\bigl(e^{i\gamma\Psi(x,y)}\sin\theta(x,y)-
e^{i\gamma\Psi_\nu(x,y)}\sin\theta_\nu(x,y)\bigr) z_j^{L,M,\nu}(x,y)
a_{j_1}(x_1+y_1),
\tag 7.5.4
\\
&Z_{\nu,j}^2(x,y)=2^{-j_1-j_2}
e^{i\gamma\Psi_\nu(x,y)}\sin\theta_\nu(x,y) z_j^{L,M,\nu}(x,y)
\bigl(a_{j_1}(x_1+y_1)-1\bigr),
\tag 7.5.5
\\
&Z_{\nu,j}^3(x,y)=2^{-j_1-j_2}
e^{i\gamma\Psi_\nu(x,y)}\sin\theta_\nu(x,y) z_j^{L,M,\nu}(x,y),
\tag 7.5.6
\endalign
$$
and form operators $\cZ_{\nu}^i =\sum_{j\in A}\cZ_{\nu,j}^i$ where
$\cZ_{\nu,j}^i$ has
 kernel $Z_{\nu,j}^i$.

The operator $\cZ_{\nu}^2$ is handled by the argument in the proof
of Proposition 5.1, with only notational changes.

The operator $\cZ_{\nu}^{3}$ represents the main term. Note
however that $\cZ_{\nu}^3 f(x) =g(x)\widetilde \cZ_{\nu,3}f(x)$,
where $g$ is a bounded function, and $\widetilde \cZ_{\nu,3}$ is an
operator which
 is already shown to be bounded
by Theorem 4.1. Thus $\|\cZ_{\nu}^3\|=O(1)$.

It remains to estimate the kernel $Z_{\nu,j}^1$. Suppose that
$a_{j_1}(x_1+y_1)z_j^{L,M,\nu}(x,y)\neq 0$. Then
$$
\align
|\Psi(x,y)-\Psi_\nu(x,y)|&\le
|p(x_1+y_1)-p(2c_\nu)-p'(2c_\nu)(x_1+y_1-2c_\nu)||x_2-y_2|
\\
&\le |x_2-y_2|\sum_{l=2}^n\frac{|p^{(l)}(2c_\nu)|}{l!}
|x_1+y_1-2c_\nu|^l
\\
&\le |x_2-y_2|\sum_{l=2}^n
\Big| \sum_{s=0}^{n-l}p^{(l+s))}(2x_1)
\frac{(2c_\nu-2x_1)^s}{s!}\Big|
\frac{|x_1+y_1-2c_\nu|^l}{l!}
\\
&\lc 2^{j_2}\sum_{l=2}^n 2^{L-(l-2)(j_1+10)} 2^{(M-L-m+1)l}\\
&\lc 2^{L+2j_1+j_2}
\sum_{l=2}^n 2^{(M-L-m+1-j_1)l}\\&\lc 2^{L+2j_1+j_2+\frac m2},
\tag 7.5.7
\endalign
$$
since we assumed that $M-L-m-j_1\le m/2n$. Similarly
$$
\align|\theta(x,y)-\theta_\nu(x,y)|&\le
\beta(x_1-y_1)|x_2-y_2||p'(2x_1)-p'(2c_\nu)|
\\
&\lc 2^{j_1+j_2}
\sum_{l=2}^n\frac{|p^{(l)}(2x_1)|}{(l-1)!}
|2x_1 -2c_\nu|^{l-1}
\\
&\lc 2^{j_1+j_2}
\sum_{l=2}^n
2^{L-(l-2)(j_1+10)} 2^{(M-L-m+1)(l-1)}
\\
&\lc 2^{L+2j_1+j_2+\frac m2} .
\tag 7.5.8
\endalign
$$

Moreover
$$
|2x_2+p'(2x_1)|\approx
|2y_2+p'(2y_1)|\approx 2^M,
\tag 7.5.9
$$
and then
$$
|2x_2+p'(2c_\nu)|\approx
|2y_2+p'(2c_\nu)|\approx 2^M
\tag 7.5.10
$$
because $j_2\ll M$ and
$$
\align
|p'(2x_1)-p'(2c_\nu)|
&\le\sum_{l=2}^n \frac{|p^{(l)}(2c_\nu)|}{(l-1)!}
|2x_1-2c_\nu|^{l-1}
\\
&\le
\sum_{l=2}^n 2^{L-j_1(l-2)}2^{(M-L-m)(l-1)}
\\
&\le 2^{L+j_1+\frac m2 +2}\le 2^{M-\frac m2 +2}.
\endalign
$$
Similarly $|p'(2y_1)-p'(2c_\nu)|\le 2^{M-\frac m2 +2}$.

{}From (7.5.7) and (7.5.8), it follows that
$$
\big|
e^{i\gamma\Psi(x,y)}\sin\theta(x,y)-
e^{i\gamma\Psi_\nu(x,y)}\sin\theta_\nu(x,y)\big|
 \lc 2^{L+2j_1+j_2+\frac m2}
$$
for the relevant values of $(x,y)$ in (7.5.4), and Schur's test
yields
$$
\|\cZ_{\nu,j}^1\|\lc 2^{L+2j_1+j_2+\frac m2}.
$$
By (7.5.9) and (7.5.10), we may apply either Proposition 5.4 with the
polynomial $p$ or with the affine polynomial
$p(c_\nu)+p'(c_\nu)(s-c_\nu)$, and the suitable choice of $\rho_j$
in part (iii) of Proposition 5.4. This leads to
$$
\|\cZ_{\nu,j}^1\|\lc \min\{ 2^{-(M+j_1+j_2)/4},
 2^{M+j_1+j_2}\}
\tag 7.5.12
$$
if $j_1\le M-L$.  We obtain
$$
\sum\Sb j\in I\\
M+j_1+j_2=\ell\endSb
\bigl\|\cZ_{\nu,j}^1\bigr\|
\lc
\sum\Sb
j_1\le M-L-m\\
L+2j_1+j_2\le 0\\
M+j_1+j_2=\ell\\
\endSb\min\{ 2^{L+2j_1+j_2+\frac{m}{2}},
 2^{-(M+j_1+j_2)/4}\}\lc \min\{ 2^{-\frac{\ell}4}, 2^{\ell-\frac m 2}\}.
$$
Summing over $\ell$ demonstrates the boundedness of the operator
$\cZ_\nu^1$.

We have shown that $\cZ^I_\nu=\sum_{i=1}^3\cZ_{\nu}^i$ is bounded
with operator norm uniformly in $M,L,m,I,\nu$.  The assertion of
the theorem now follows from (7.5.1).\qed

\head{\bf 8. Failure of weak amenability for Lie groups}
\endhead

Suppose that $G$ is a connected Lie group, with Lie algebra $\mfr{g}$.
Then $\mfr{g}$ decomposes as $\mfr{s}\oplus\mfr{r}$, where $\mfr{s}$
is a semisimple subalgebra and $\mfr{r}$ is the maximal solvable ideal
of $\mfr{g}$. We may write $\mfr{s}$ as a sum of simple ideals:
$$
\mfr{s}=\mfr{s}_1\oplus\cdots\oplus\mfr{s}_m.
\tag 8.1
$$
Denote by $R$, $S$ and $S_i$ the analytic subgroups of $G$ corresponding
to $\mfr{r}$, $\mfr{s}$ and $\mfr{s}_i$. Then $R$ is closed, but $S$ and $S_i$
need not be. Further, $G=SR$, but this need not be a semidirect
product, as $S\cap R$ may be nontrivial. To do analysis on $G$, we need
$S$ to be closed and the product $SR$ to be semidirect.  Our first result
enables us to work in this better environment by passing
to a finite covering group.

\proclaim{Proposition 8.1}
Let $G$, $R$, $S$ and $S_i$ be as described above, and suppose that
$S$ has finite center.  
 Then $G$ has a finite  covering group
$G^\natural$ which has
 closed connected subgroups $R^\natural$, $S^\natural$ and
$S_i^\natural$, whose Lie algebras are
$\mfr{r}$, $\mfr{s}$ and $\mfr{s}_i$,
such that $R^\natural$ is normal and solvable, $S^\natural$ is the direct
product of the  simple Lie groups  $S_i^\natural$, each of which has finite
center, and $R^\natural\cap S^\natural=\{e\}$; thus $G^\natural$ is the
semidirect product $S^\natural\ltimes R^\natural$.
\endproclaim

Here $G^\natural$ is said to be 
 a finite covering group if $G$ is isomorphic to $G^\natural/Z$ where
$Z$ is a finite normal subgroup of $G^\natural$.

We leave the proof of Proposition 8.1 until later.
Observe that
$\Lambda(G)=\Lambda(G^\natural)$, by 1.2.4, (ii); moreover by 1.2.4,
(i),
$G^\natural$ has
a multiplier bounded approximate unit
if $G$ has one.
Thus to compute $\Lambda(G)$, we
may and shall henceforth assume that $G,R,S$ and $S_i$ have the properties of
$G^\natural,R^\natural,S^\natural$ and $S_i^\natural$ in Proposition 8.1.

Now to prove the theorem, observe that if the factors $S_i$
making up $S$ are all either compact (when $i\in I$, say)
 or of real rank one and commute with $R$
(when $i\in J$, say) then we may write $G$ as a direct product:
$$
G=(\prod_{i\in J}S_i)\times\bigl((\prod_{i\in I}S_i)
\ltimes R\bigr).
\tag 8.2
$$
The second factor is amenable, so $\Lambda((\prod_{i\in
I}S_i)\ltimes R)=1$, and hence
$$
\Lambda(G)=\prod_{i\in J}\Lambda(S_i)=\prod_{i=1}^m\Lambda(S_i),
\tag 8.3
$$
by 1.2.1(iv).

On the other hand, if any $S_i$, $i\in J$,  is of
real rank at least two, then
$\Lambda(G)\geq\Lambda(S_i)=+\infty$ (see
 1.2.1(iii));
moreover the proof in \cite{14} and  \cite{10}
that $\Lambda (S_i)=\infty$
in combination with (1.2.4 (i))
  shows that
$S_i$ and therefore $G$ does not have multiplier bounded
 approximate units.

The remaining case to consider
is when there is a factor $S_i$ of real rank one
which does not centralize $R$. The following result contains
the structural
information needed to reduce to known cases.

\proclaim{Proposition 8.2} Suppose that the connected Lie group $G$ is
a semidirect product of the form $S\ltimes R$, where $S$ is closed,
connected, semisimple and has finite center, and $R$ is closed, connected
and solvable, and suppose that a noncompact factor $S_i$ of $S$ does
not centralize $R$.   Then $G$ contains a closed subgroup $G_0$ with a
compact normal subgroup $K_0$ such that $G_0/K_0$, or a double cover
of $G_0/K_0$, is isomorphic to $\SLtwoR\ltimes\R^n$ (where $n\geq 2$) or
to $\SLtwoR\ltimes\Heis^n$ (where $n\geq 1$).
\endproclaim

Thus under the assumptions of Proposition 8.2. it follows
 that $G_0/K_0$ does not admit multiplier
bounded approximate units,
by the calculations of \cite{9} for the groups $\SLtwoR\ltimes\R^n$
and of
this paper for the groups $\SLtwoR\ltimes\Heis^n$. Thus by (1.2.1), (iii)
and (1.2.4), (i) the group  $G$ does not
have multiplier
bounded approximate units and in particular we have  $\Lambda(G)=\infty$.

It remains to prove
Propositions 8.1 and 8.2.

\demo{\bf Proof of Proposition 8.1}   To every Lie algebra $\mfr{a}$, we may
associate a unique connected, simply connected Lie  group $A$ with Lie
algebra $\mfr{a}$.   Every connected Lie group $A'$ with Lie algebra
$\mfr{a}$ is a quotient of $A$ by a discrete normal, and hence central,
subgroup $D$ of $A$.  For these facts, and much more, about the
structure of Lie algebras and Lie groups, see, for instance, \cite{11} or
\cite{21}, \cite{34}.  Consequently, we will be
interested in the structure of the
center of a connected, simply connected Lie group.

Let $G^\sharp$ be the simply connected covering group of $G$, and let
$R^\sharp$ and $S^\sharp$ be the subgroups of $G^\sharp$
corresponding to $\mfr{r}$ and $\mfr{s}$.
Then $R^\sharp$ and $S^\sharp$ are both closed in $G^\sharp$; further,
$R^\sharp$ is normal in $G^\sharp$ and $S^\sharp\cap R^\sharp=\{e\}$,
so that $G^\sharp=S^\sharp\ltimes R^\sharp$ (\cite{34}, Thm.~3.18.13).
Consequently, the center $Z(G^\sharp)$ of $G^\sharp$ may be written
as a direct product: $Z(G^\sharp)=Z_S\times Z_R$, where $Z_R$ is the
subgroup of the center $Z(R^\sharp)$ of $R^\sharp$ of elements which
commute with $S^\sharp$, and $Z_S$ is the subgroup of the center
$Z(S^\sharp)$ of $S^\sharp$ of elements commuting with $R^\sharp$.
Let $S_i^\sharp$ be the subgroup of $G^\sharp$ corresponding to $\mfr{s}_i$.
Then $S_i^\sharp$ is closed in $S^\sharp$, and hence in $G^\sharp$; further,
$S_i^\sharp$ is simply connected and normal in $S^\sharp$, so that  $S^\sharp$
is a direct product of the factors $S_i^\sharp$, and $Z(S^\sharp)$ is the
direct product of the centers $Z(S_i^\sharp)$  of the $S_i^\sharp$
(\cite{34}, Thm.~3.18.1).

The group $G$, being a quotient of $G^\sharp$, is of the form $G^\sharp/D$,
where $D$ is a discrete subgroup of $Z_S\times Z_R$.  Set
$$D_0=\prod_{i=1}^m [D\cap S_i^\sharp] \times [D\cap R^\sharp].$$
We need an auxiliary result.
\proclaim {Lemma 8.1.1}

(i)  $D_0$ is of finite index in $D$.

(ii) Each
$D\cap S_i^\sharp$ is of finite index in the center $Z(S_i^\sharp)$
of $S_i^\sharp$.
\endproclaim

Taking this for granted  the group $G/D_0$ is a finite covering  of $G/D$,
and has  the required properties as it is  isomorphic to
$$
\prod_{i=1}^m [S_i^\sharp/ (D\cap S_i^\sharp) ] \ltimes [R^\sharp/(D\cap R^\sharp)].
$$
We take $S_i^\natural$ to be $S_i^\sharp/ (D\cap S_i^\sharp)$ and $R^\natural$
to be $R^\sharp/(D\cap R^\sharp)$.  Then $R^\natural$ is closed, normal
and solvable, and $S_i^\natural$ is closed, simple and have finite center.
The center of $S^\natural$ is the product of the centers of the
groups $S_i^\natural$, and is also be finite.  Finally, $S^\natural\cap R^\natural$ is trivial, where $S^\natural = \prod_{i=1}^m S_i^\natural$.

It remains to give the
\demo{Proof of Lemma 8.1.1}
First, we claim that  $Z_S$ is of finite index in $Z(S^\sharp)$.  Indeed,
the adjoint action $\Ad_{\mfr{r}}$ of $S^\sharp$ on  $\mfr{r}$ is a linear
representation of $S^\sharp$, and the image $\Ad_{\mfr{r}}S^\sharp$
of $S^\sharp$ in $\text{SL}(\mfr{r})$ is a closed semisimple subgroup
of $\text{SL}(\mfr{r})$;  the center $C$ of this subgroup is finite
\cite{21, Prop.~7.9}  Moreover, $\Ad_{\mfr{r}}(Z(S^\sharp))$
is contained in $C$, so that  $\Ad_{\mfr{r}}(Z(S^\sharp))$ is finite.
The group $R^\sharp$ is generated by arbitrarily small neighborhoods
of the identity, so an element of $S^\sharp$ centralizes $R^\sharp$ if and
only if it centralizes small neighborhoods of the identity, and hence if and
only if it acts trivially on $\mfr{r}$ by the adjoint action.  Then $Z_S$ is
the kernel of the adjoint map of $Z(S^\sharp)$ into $\text{SL}(\mfr{r})$.
In conclusion, $\Ad_{\mfr{r}}(Z(S^\sharp))$ is isomorphic to $Z(S^\sharp)/Z_S$,
so $Z_S$ is indeed of finite index in $Z(S^\sharp)$.

Next, note that $D\cap S^\sharp = D\cap Z(S^\sharp) = D\cap Z_S$, since
$D \subseteq Z_S \times Z_R$.  We claim that $D\cap S^\sharp$ is of finite
index in $Z(S^\sharp)$.   Indeed, $S$ is isomorphic to $S^\sharp/D\cap S^\sharp$, and its center is isomorphic to $Z(S^\sharp)/D\cap S^\sharp$ (in fact, if
$x\in S^\sharp \setminus Z(S^\sharp)$, then $\Ad(x)$ acts nontrivially on
any small neighborhood of the identity in $S^\sharp$, and hence its image in
$S^\sharp/D\cap S^\sharp$ also acts nontrivially on such a neighborhood).
By hypothesis, the center of $S$ is finite, so $D\cap S^\sharp$ is of
finite index in $Z(S^\sharp)$.  Since
$$Z(S_i^\sharp)/[D\cap S_i^\sharp] \simeq DZ(S_i^\sharp)/D
\subseteq DZ(S^\sharp)/D \simeq Z(S^\sharp)/[D\cap S^\sharp],
$$
$D\cap S_i^\sharp$ is of finite index in $Z(S_i^\sharp)$.
Now write $E$ for $\prod_{i=1}^m [D\cap S_i^\sharp]$; it follows that
$E$ is a subgroup of $Z(S^\sharp)$ of finite index which is part (ii) of
the lemma.

For part (i) we shall establish that  $E (D\cap R^\sharp)$ is of finite index in $D$.  This also
follows from the isomorphism theorems.  Since $Z_S/(D\cap Z_S)
\subset Z(S^\sharp)/(D\cap Z_S)$ is finite, and
$ Z_S/(D\cap Z_S) \simeq D Z_S/D$, there exist $z_1,z_2,\ldots,z_J\in Z_S$
such that
$$
DZ_S = \bigcup_{j=1}^J Dz_j.
$$
If $Dz_j\cap Z_R \neq \emptyset$, take $r_j \in Dz_j\cap Z_R$; otherwise take
$r_j = e$.  Then it is easy to check that
$$
Dz_j \cap Z_R \subseteq (D\cap Z_R) r_j .
$$
It follows that
$$
(DZ_S)\cap Z_R
= \bigcup_{j=1}^J Dz_j \cap Z_R
\subseteq  \bigcup_{j=1}^J (D \cap Z_R)r_j,
$$
and so $(DZ_S)\cap Z_R \big/(D\cap Z_R)$ is also finite.  As $D$ is contained in the direct product $Z_S Z_R$, it follows that $DZ_S = Z_S((DZ_S)\cap Z_R)$,
and so
$$\align
D\Big/(D\cap Z_S)(D\cap Z_R)
&\subseteq DZ_S \Big/ (D\cap Z_S)(D\cap Z_R)  \\
&= Z_S((DZ_S)\cap Z_R)\big/ (D\cap Z_S)(D\cap Z_R)  \\
&\simeq (Z_S/(D\cap Z_S)) ((DZ_S)\cap Z_R/D\cap Z_R),  \\
\endalign$$
which is finite.  Thus $(D\cap Z_S)(D\cap Z_R)$ is of finite index in $D$.
Since $E$ has finite index in $D\cap Z_S$, $E(D\cap Z_R)$ is of finite index
in $(D\cap Z_S)(D\cap Z_R)$, and thus of finite index in $D$. This finishes the proof of Lemma 8.1.1 and Proposition  8.1.
\qed\enddemo
\enddemo

\demo{\bf Proof of Proposition 8.2}

We recall that $S_i$ does not centralize $R$ so that $\frak s_i$
does not centralize $\frak r$.  The proof now  proceeds by a series of reductions.

First, we cut down the
semisimple part. Take a Cartan involution $\theta$ of $\mfr{s}_i$,
 so that
$\mfr{s}_i = \mfr{k}\oplus\mfr{p}$, where $\mfr{k}$ and $\mfr{p}$ are the
$+1$ and $-1$ eigenspaces of $\theta$, and take a maximal abelian
subalgebra $\mfr{a}$ of $\mfr{p}$.  Then the Lie algebra $\mfr{s}_i$
decomposes into a sum:
$$
\mfr{s}_i = \mfr{g}_0 + \sum_{\alpha \in \Sigma} \mfr{g}_\alpha ,
$$
where each $\alpha$ is a linear functional on $\mfr{a}$, and
$X\in \mfr{g}_\alpha$ if and only if $[H,X] = \alpha(H)X$ for all $H\in\mfr{a}$.
For more on these root  decompositions, see, e.g., \cite{15},
\cite{21}.
Take
a nonzero element $X$ in this $\mfr{g}_\alpha$,
where $\alpha \neq 0$.
Then $\spaN\{X,\theta X,[X,\theta X]\}$
is a subalgebra $\mfr{s}_0$ of $\mfr{s}_i$ isomorphic to $\mfr{sl}(2,\R)$, and
the corresponding analytic subgroup $S_0$ of $S_i$ is locally
isomorphic to
$\SLtwoR$, has finite center, does not centralize $R$, and is closed in $S_i$ and hence in $S$ (see
\cite{35, Lemma~1.1.5.7}).

Note also that $S_0/Z(S_0)$ is isomorphic to a matrix group (\cite{34, Thm.~2.13.2}) and the only matrix groups locally isomorphic are
$SL(2,\Bbb R)$
and $PSL(2,\Bbb R) $
(i.e., $SL(2,\Bbb R)$ divided by its center).

  Furthermore, $[\mfr{s}_0,\mfr{r}]\ne 0$; indeed,
$\{X\in\mfr{s}_i:\ad(X)|_{\mfr{r}}=0\}$
is an ideal in
$\mfr{s}_i$, which is a simple Lie algebra;  and hence
$\{X\in\mfr{s}_i:\ad(X)|_{\mfr{r}}=0\}=\{0\}$. The subgroup $S_0\ltimes R$
is closed in $G$.

The second reduction cuts down to the nilradical. Let $N$ be the maximal
connected normal nilpotent subgroup of $G$, which is automatically closed, and
let $\mfr{n}$ be its Lie algebra. Then $N\subseteq R$ and $\mfr{n}\subseteq
\mfr{r}$, and moreover, $[\mfr{s}_0,\mfr{r}]\subseteq\mfr{n}$ (see
\cite{34, Thm.~3.8.3}).
We claim that
$[\mfr{s}_0,\mfr{n}]\neq \{0\}$.
If it were true that $[\mfr{s}_0,\mfr{n}]=\{0\}$, then the Jacobi
identity would imply that,
$$
[[X,Y],Z]=[[X,Z],Y]+[X,[Y,Z]]=0,\quad
\text{  $X, Y\in \mfr{s}_0$,  $Z\in \mfr{r}$,}
$$
since the inner commutators of both summands
of the middle term of the equality lie in $\mfr{n}$, from which it would follow
that $[\mfr{s}_0,\mfr{r}]=\{0\}$. Thus $[\mfr{s}_0,\mfr{n}]\ne\{0\}$. It now suffices to
consider $S_0\ltimes N$, which is closed in $S_0\ltimes R$ and hence in $G$.


The third reduction allows us to assume that $N$ is simply connected.
Let $K$ be the maximal compact connected central subgroup
of $S_0\ltimes N$; it is
contained in the nilradical $N$. We observe that the nilradical of
$(S_0\ltimes N)/K$ is equal to
$S_0\ltimes (N/K)$ and we shall show that $N':=N/K$ is simply connected.


The center $Z(N)$ of $N$  is a connected Abelian Lie group
 \cite{34, Cor.~3.6.4}
and thus isomorphic to $\Bbb R^k\times\Bbb T^l$,
for suitable $k,l$.
We claim that
 $\Bbb T^l$ is central in $S_0\ltimes N$.
Note that then $\bbT^l$ is also the maximal compact connected central subgroup
of $G$ (since any connected central subgroup must be in the nilradical).

We first show that $\bbT^l$ is a normal subgroup.
 For any fixed $g\in S_0\ltimes N$ the automorphism
$\phi_g: n\mapsto gng^{-1}$ fixes the center; thus $g\bbT^l g^{-1}$ is a
compact Lie  subgroup of $Z(N)$ which must be $\Bbb T^l$; thus $\bbT^l$ is normal in $G$.
Consider the map $g\to \phi_g|_{\bbT^l}$ which takes $G$ into the automorphism
group of $\Bbb T^l$. This group is discrete and since $G$ is connected we see that
$\phi_g|_{\bbT^l}$ is the identity, thus $\Bbb T^l$ is central in $G$
and thus isomorphic to $K$.

We claim that $N':=N/K$ is simply connected. Indeed let $\widetilde N$ be the simply connected covering group of $N$; it has center $Z(\widetilde N)=\Bbb R^{k+l} \supset \Bbb Z^l$ and $N=\widetilde N/\Bbb Z^l$. Now
$$
N'=N/\Bbb T^l= (\widetilde N/\Bbb Z^l)/(\Bbb R^l/\Bbb Z^l)\cong \widetilde N/\Bbb R^l.
$$
But $N'\cong \widetilde N/\Bbb R^l$ is simply connected by
\cite{34, Thm.~3.18.2}.

In the remaining part of the proof we shall show that
$S_0\ltimes N'\cong (S_0\ltimes N)/K$ has
a closed subgroup $G_1$ which is locally isomorphic to
some
$\SLtwoR\ltimes\R^n\;(n\geq 2)$ or some
 $\SLtwoR\ltimes\Heis^n\;(n\geq1)$.
The
desired subgroup $G_0$ of $G$ is then the closed subgroup of
$S_0\ltimes N$ of
elements whose image in $(S_0\ltimes N)/K$ under the canonical
projection lies
in $G_1$, and the appropriate compact subgroup  $K_0$
is the direct product of $Z(S_0\ltimes N)\cap
S_0$ and $K$.

In order to proceed we need the following lemma.

\proclaim{Lemma 8.2.1}
Let $\pi_n:\mfr{sl}(2,\R)\to\End(\R^n)$ be the (unique)
irreducible representation of $\mfr{sl}(2,\R)$ of dimension $n$.
The space of bilinear forms $B:\R^n\times\R^n\to\R$ which satisfy
$$
B(\pi_n(U)V,W)+B(V,\pi_n(U)W)=0
\qquad\forall U\in\mfr{sl}(2,\R)\quad\forall V,W\in\R^n  \tag 8.4
$$
is one-dimensional. These forms are symmetric or skew-symmetric
as $n$ is odd or even.
\endproclaim

\demo{Proof}
Let
$$\text{
$H=\pmatrix 1&0\\0&-1\endpmatrix$,
$X=\pmatrix 0&1\\0&0\endpmatrix$,
$Y=\pmatrix 0&0\\1&0\endpmatrix$}$$
that is, the standard basis for $\mfr{sl}(2,\R)$
satisfying the commutation relation
$[X,Y]=H$, $[H,X]=2X$ and $[H,Y]=-2Y$. It is well known (see for instance
\cite{18, ch.~III.8})
 that there is a basis
$\{E_0,\dots,E_{n-1}\}$ for $\R^n$
such that
$$\aligned
\pi_n(H)E_j&=(n-1-2j) E_{j}, \quad j=0,\dots, n-1
\\
\pi_n(X) E_j&=E_{j+1},
 \quad j=0\dots, n-2, \quad \pi_n(X)E_{n-1}=0
\\
\pi_n(Y) E_j&=j(n-j+1)E_{j-1},
 \quad j=1\dots, n-1, \quad \pi_n(Y)E_{0}=0.
\endaligned
$$

Let $B$ be a bilinear form satisfying (8.4). If $0\leq i,j\leq n-1$, then
$$
0= B(\pi_n(H)E_i,E_j)+B(E_i,\pi_n(H)E_j)= (2n-2i-2j-2)B(E_i,E_j)
$$
so that
$$B(E_i,E_j)=0 \text{ if }j\neq n-i-1.
\tag 8.5
$$ Further, if $1\le j\le n-1$, then
$$
  B(E_{j},E_{n-j-1})
= B(\pi_n(X)E_{j-1},E_{n-j-1})
=-B(E_{j-1},\pi_n(X)E_{n-j-1})
=-B(E_{j-1},E_{n-j}),
$$
whence
$$
B(E_i,E_{n-i-1})=(-1)^{j}B(E_0,E_{n-1}), \qquad
1\leq i\leq n-1,
$$
 so that $B$ is completely determined by $B(E_0,E_{n-1})$. In
particular,
$$B(E_i,E_{n-i-1})=(-1)^{n-1}B(E_{n-i-1},E_i).
\tag 8.6$$
Thus by (8.5) and (8.6),  $B$ is symmetric if $n$ is odd and skew-symmetric if $n$ is even.
\qed\enddemo

\demo{Proof of Proposition 8.2, continued}
We now consider $S_0\ltimes N'$
 and we must produce a closed subgroup of $S_0\ltimes N'$
locally isomorphic to $\SLtwoR\ltimes\R^n\;(n\geq 2)$ or to
$\SLtwoR\ltimes\Heis^n\;(n\geq1)$. Let $\frak n$ be the Lie algebra of
$N'$; since $N'$ is simply connected, the
exponential map is a homeomorphism from $\mfr{n}$ to $N'$, and subalgebras of
$\mfr{n}$ map to closed subgroups of $N'$.

We define the ascending central series of $\mfr{n}$ inductively: let
$\mfr{n}_0$ be $\{0\}$, and if $j\geq 1$, define $\mfr{n}_j$ to be
$\{X\in\mfr{n}:[X,\mfr{n}]\subseteq \mfr{n}_{j-1}\}$.
Since $\mfr{n}$ is nilpotent, there exists a positive integer $l$
such that $\mfr{n}_l=\mfr{n}$, so
$$
\{0\}=\mfr{n}_0\subset\mfr{n}_1\subset\dots\subset\mfr{n}_l=\mfr{n}.
$$
Choose $j$ such that $[\mfr{s}_0,\mfr{n}_{j-1}]=\{0\}$ but
$[\mfr{s}_0,\mfr{n}_j]\ne\{0\}$. Under the action of the semisimple group $S_0$
on $\mfr{n}$, the subalgebra $\mfr{n}_j$ splits into a sum of irreducible
$\Ad(S_0)$ modules, not all of which are trivial. Let $\mfr{m}$ be a
nontrivial summand in this decomposition; then $[\mfr{s}_0,\mfr{m}]=\mfr{m}$.

{}From the Jacobi
identity   we get
$$[[X,Y],Z]=[[X,Z],Y]+[X,[Y,Z]]=0,\quad
\text{  $X\in\mfr{s}_0,\;  Y\in\mfr{m}$ and $Z\in\mfr{n}_{j-1}$},
$$
since $[X,Z]\in[\mfr{s}_0, \mfr{n}_{j-1}]=\{0\}$ and
$[X,[Y,Z]]\in[\mfr{s}_0,\mfr{n}_{j-1}]=\{0\}$. It follows that
$[\mfr{m,n}_{j-1}]=0$. In particular, $[\mfr{m,m}]\subseteq\mfr{n}_{j-1}$, so
$[\mfr{m,[m,m}]]=\{0\}$, and $\mfr{m+[m,m}]$ is a subalgebra of $\mfr{n}$.
Given any linear form $\lambda$ on $\mfr{[m,m]}$, the bilinear form
$B:(V,W)\mapsto\lambda[V,W]$ satisfies
$$
 B(\ad(U)V,W)+B(V,\ad(U)W)
=\lambda([\ad(U)V,W]+[V,\ad(U)W])
=\lambda(\ad(U)[V,W])
=0
$$
for all $U$ in $\mfr{sl}(2,\R)$ and all $V$ and $W$ in $\mfr{m}$. Since the
space of such bilinear forms is one-dimensional, from Lemma 8.2.1, it follows
that $\dim([\mfr{m,m}])\leq 1$ (in particular, if $\dim(\mfr{m})$ is odd, then
$[\mfr{m}, \mfr{m}] = \{0\}$, for the form $(V,W)\mapsto\lambda[V,W]$ is
skew-symmetric). Let $m=\text{dim}( \frak m$).
Then $\exp(\mfr{m+[m,m}])$ is isomorphic to $\Bbb R^m$ if
$m$ is odd or $m$ is even and $[\frak m,\frak m]=0$, and isomorphic
to $H^{m/2}$ if $m$ is even and $[\frak m,\frak m]\neq 0$.
The group
$S_0\ltimes\exp(\mfr{m+[m,m}])$ is the required subgroup of
$S_0\ltimes N'$.
\qed\enddemo
\enddemo

\enddocument